\documentclass[12pt]{article}
\usepackage{amsmath,amssymb}
\usepackage{amsfonts}
\usepackage{multicol} 
\usepackage{array} 
\usepackage{lscape} 
\usepackage{xcolor} 
\usepackage{enumerate}
\usepackage{pictex}
\usepackage{appendix}

\newtheorem{Def}{Definition}[section]

\newtheorem{Lema}[Def]{\bf Lemma}

\newtheorem{Teo}[Def]{\bf Theorem}

\newcommand{\too}{\rightarrowtail}
\newcommand{\landb}{\land_\text{\tiny B}}
\newcommand{\lorb}{\lor_\text{\tiny B}}
\newcommand{\simb}{\sim_\text{\tiny B}}

\newcommand{\Br}{\mbox{$\cal B$}} 

\newenvironment{proof}{\noindent\bf Proof. \rm}{\hfill$\blacksquare$}

\newenvironment{proofish}{\noindent\bf Proof. \rm}{}

\newcommand{\gen}{\reflectbox{$\neg$}}

\renewcommand{\lceil}{\gen}

\title{An alternative axiomatic presentation of Nelson algebras}
\author{Juan Manuel Cornejo \and Andr\'es Gallardo \and Luiz Monteiro \and Ignacio Viglizzo}

\date{}
\begin{document}

\maketitle

\begin{abstract}
Nelson algebras are defined in \cite{brignole67caracterisation} in terms of the language $\langle \land,\lor,\to,\sim,1\rangle$.
In 1962, Diana Brignole \cite{Brignole65Nlattice}, solved the problem proposed by Antonio Monteiro, of giving an axiomatization of Nelson algebras in terms of the connectives $\too$, $\land$ and the constant $0=\sim 1$. In this work we present for the first time a complete proof of this fact, and also show the dependence and independence of some of the axioms proposed by Brignole.
\end{abstract}

\section{Preliminaries}

Nelson algebras or $\mathcal{N}$-algebras were introduced by H.~Rasiowa \cite{rasiowa58nlattices} as an algebraic counterpart of Nelson's constructive logic with strong negation \cite{nelson1949constructible}. Later, D.~Brignole and A.~Monteiro \cite{brignole67caracterisation}, \cite{brignole69equational}, gave a characterization using identities, proving that Nelson algebras form a variety. This characterization was given in terms of the operations $\land,\lor,\to,\sim$ and the constant $1$. 

A different implication operation can be defined by: 
\[x\too y=(x\to y)\land(\sim y\to \sim x)\]

In 1962, Diana Brignole, solved the problem proposed by A.~Monteiro, of giving an axiomatization of Nelson algebras in terms of the connectives $\too$, $\land$ and the constant 0. This solution was communicated in the annual meeting of the Uni\'on Matem\'atica Argentina, and a summary (with some typos) was published in \cite{Brignole65Nlattice}, but the corresponding proof, to our knowledge, has not been published.

In \cite{spinks07characterisations}, Spinks and Veroff  used this axiomatization to prove that the variety of Nelson algebras is term equivalent to a variety of bounded 3-potent BCK-semilattices, and in \cite{spinks08constructiveI} and \cite{spinks08constructiveII},  they prove that using the operation $\too$, Nelson algebras can be understood as residuated lattices, with the product given by the term:
\[x*y=\sim(x \to\sim y) \lor\sim(y \to\sim x). \] As a consequence, the corresponding logic, constructive logic with strong negation can be seen as a substructural logic (see also \cite{busaniche10constructive}).

In this note we give a complete proof of the axiomatization proposed by Brignole, prove the independence of some of the axioms from the rest and also announce that two of the identities can be derived from the others, although we only have an automated proof of this fact.

\begin{Def}\label{defnel}

A \emph{Nelson algebra} is an algebra $\langle A,\land,\lor,\to,\sim,1 \rangle$ of type $(2,2,2,1,0)$ such that the following conditions are satisfied for all $x,y,z$ in $A$:
	
\begin{enumerate}[\rm (N1)]
	\item $x\land(x\lor y)=x$,\label{primercondicretN}
	\item $x\land(y\lor z)=(z\land x)\lor(y\land x)$,\label{segundacondicretN}
	\item $\sim\sim x=x$,\label{involucionN}
	\item $\sim(x\land y)=\sim x\lor \sim y$,\label{deMenInfN}
	\item $x\land \sim x=(x\land \sim x)\land (y\lor \sim y)$,\label{kleeneN}
	\item $x\to x=1$,\label{ximpxigual1N}
	\item $x\land(x\to y)=x\land(\sim x\lor y)$,\label{relacImplicaInfN}
	\item $(x\land y)\to z=x\to(y\to z)$.\label{relacImplicaInf2N}
\end{enumerate}
\end{Def}

We will denote by $\mathcal{N}$ the variety of Nelson algebras.

The axioms in this list form an independent set, see \cite{monteiro96axiomes}.

By axioms (N\ref{primercondicretN}) and (N\ref{segundacondicretN}) we have that every Nelson algebra is a distributive lattice (see Sholander \cite{sholander51postulates}). Furthermore, if we define $0=\sim 1$, we have that 0 and 1 are the bottom and top element of $A$ respectively.

In a Nelson algebra we can also define the following operations that will be used in this work:

\begin{itemize}
	\item $\lceil x:=x\to(\sim 1)$.
	\item $x\too y:=(x\to y)\land(\sim y\to\sim x)$.
\end{itemize}

\begin{Lema} \label{lema_prop_Nelson}
Let $\langle A,\land,\lor,\to,\sim,1\rangle$ be a Nelson algebra. The following properties are satisfied in $A$:
\begin{enumerate}[\rm (a)]
	\item $x\to(y\land z)=(x\to y)\land(x\to z)$, \label{distribInfADerN}
	\item $1\to x=x$,\label{tesis_1_6}
	\item $\sim x\leq\lceil x$,\label{tesis_1_61}
	\item $(x\to x)\land(\sim x\to\sim x)=1$,\label{tesis_1_29}
	\item $\sim y\leq y\to z$,\label{tesis_1_8}
	\item $y\leq x\to y$,\label{tesis_1_2}
	\item $(x\lor y)\to z=(x\to z)\land(y\to z)$,\label{tesis_1_9}
	\item $x\to y=x\too(x\too y)$,\label{tesis_1_83}
	\item $x\to(x\to y)=x\to y$,\label{tesis_1_34} 
\end{enumerate}
\end{Lema}

\begin{proof} 
The proofs of these items  can be found in \cite{viglizzo99algebras}.
\end{proof}

\

The following definition is based on the set of equations  given by Brignole in \cite{Brignole65Nlattice}, with some typos corrected.

\begin{Def}
A \emph{Brignole algebra} is an algebra ${\bf A}=\langle A,\landb,\too,0\rangle$ of type $(2,2,0)$, such that the following equations are satisfied for all $x,y,z\in A$:

	\begin{enumerate}[\rm (B1)]\label{defbrig}
		\item $(x\too x)\too y=y$,\label{AxiomaBrignole01}
		\item $(x\too y)\landb y=y$,\label{AxiomaBrignole02}
		\item $x\landb\simb(x\landb\simb y)=x\landb(x\too y)$,\label{AxiomaBrignole03}
		\item $x\too(y\landb z)=(x\too y)\landb(x\too z)$,\label{AxiomaBrignole04}
		\item $x\too y=\simb y\too\simb x$,\label{AxiomaBrignole05}
		\item $x\too(x\too(y\too(y\too z)))=(x\landb y)\too((x\landb y)\too z)$,\label{AxiomaBrignole06}
		\item $\simb(\simb x\landb y)\too(x\too y)=x\too y$,\label{AxiomaBrignole07}
		\item $x\landb(x\lorb y)=x$,\label{AxiomaBrignole08}
		\item $x\landb(y\lorb z)=(z\landb x)\lorb(y\landb x)$,\label{AxiomaBrignole09}
		\item $(x\landb\simb x)\landb (y\lorb\simb y)=x\landb\simb x$.\label{AxiomaBrignole10}
	\end{enumerate}
	where
	\begin{itemize}
		\item $\simb x:=x\too 0$.
		\item $x\lorb y:=((x\too 0)\landb(y\too 0))\too 0$.
	\end{itemize}	
\end{Def}

We will denote by $\mathcal{B}$ the variety of Brignole algebras.

We are going to show that $\mathcal{B}$ and $\mathcal{N}$ term equivalent.

\section{Term equivalence between $\mathcal{B}$ and $\mathcal{N}$}

Let us consider a Nelson algebra $\mathbf A=\langle A,\land,\lor,\to,\sim,1\rangle$. Over $\mathbf A$ we can define the terms

\begin{itemize}
	\item $x\too y:=(x\to y)\land(\sim y\to \sim x)$,
	\item $x\landb y:=x\land y$,
	\item $0:=\sim 1$,
	\item $\simb x:=x\too 0$,
	\item $x\lorb y:=((x\too 0)\landb (y\too 0))\too 0$.
\end{itemize}

We are going to prove that $\langle A,\landb,\too,0 \rangle$ is a Brignole algebra. In order to see that, we need the following result:

\begin{Lema} \label{lema_prop_Nelson_Nuevas}
In a Nelson algebra $\mathbf A$ the following identities hold for all $x,y,z\in A$:
\begin{enumerate}[\rm (a)]
	\item $\simb x=\sim x$, \label{proposicion1} 
	\item $x\lorb y=x\lor y$, \label{proposicion2}
	\item $x\land(x\lorb y)=x$, 
	$x\land(y\lorb z)=(z\land x)\lorb(y\land x)$,
	$(x\land\simb x)\land (y\lorb\simb y)=x\land\simb x$,\label{axiomasReticuladoB}
	\item $x\too x=1$,\label{lema2_1}
	\item $x=x\land(\sim x\to y)$,\label{proposicion2_3}
	\item $1\too x=x$,\label{lema2_2}
	\item $(x\too x)\too y=y$,\label{axioma_A1}
	\item $(x\too y)\land y=y$,\label{axioma_A2}
	\item $x\land\sim(x\land\sim y)=x\land(x\too y)$,\label{axioma_A3}
	\item $x\too(y\land z)=(x\too y)\land(x\too z)$,\label{axioma_A4}
	\item $x\too y=\simb y\too\simb x$,\label{axioma_A5}
	\item $x\too(x\too(y\too(y\too z)))=(x\land y)\too((x\land y)\too z)$,\label{axioma_A6}
	\item $\sim(\sim x\land y)\too(x\too y)=x\too y$.\label{axioma_A7}
\end{enumerate}
\end{Lema}

\begin{proof}
\begin{itemize}
\item[(\ref{proposicion1})] 
$\simb x=x\too 0=x\too(\sim 1)=(x\to(\sim 1))\land(\sim \sim 1\to\sim x)\underset{\rm(N\ref{involucionN})}{=}(x\to (\sim 1))\land(1\to(\sim x))$ $\underset{\rm\ref{lema_prop_Nelson} (\ref{tesis_1_6})}{=}$ $(x\to(\sim 1))\land\sim x\underset{\rm\text{def.}\ \lceil }{=}\lceil x\land\sim x\underset{\rm\ref{lema_prop_Nelson} (\ref{tesis_1_61})}{=}\sim x$.

\item[(\ref{proposicion2})] 
$x\lorb y=((x\too 0)\land(y\too 0))\too 0\underset{\text{def.}\ \sim }{=}\sim(\sim x \land \sim y)\underset{\rm(N\ref{deMenInfN})}{=}\sim\sim x\lor\sim\sim y\underset{\rm(N\ref{involucionN})}{=}x\lor y$.

\item[(\ref{axiomasReticuladoB})] From item (\ref{proposicion2}) and the definition of $\landb$ we have that this result is immediate from axioms (N\ref{primercondicretN}), (N\ref{segundacondicretN}) and (N\ref{kleeneN}).

From the first two conditions of this item we conclude that $\langle A;\landb,\lorb\rangle$ is a distributive lattice \cite{sholander51postulates}. Therefore, for the rest of the proof we will use of properties of distributive lattices without making explicit mention them.

\item[(\ref{lema2_1})] It follows immediately from Lemma \ref{lema_prop_Nelson} (\ref{tesis_1_29}).

\item[(\ref{proposicion2_3})] It is a consequence of Lemma \ref{lema_prop_Nelson} (\ref{tesis_1_8}) and (N\ref{involucionN}).

\item[(\ref{lema2_2})] $(1\to y)\land(\sim y\to \sim1)\underset{\rm\ref{lema_prop_Nelson}(\ref{tesis_1_6})}{=}y\land(\sim y\to \sim1)\underset{\rm(\ref{proposicion2_3})}{=}y.$

\item[(\ref{axioma_A1})] $(y\too y)\too x\underset{\rm(\ref{lema2_1})}{=}1\too x\underset{\rm(\ref{lema2_2})}{=}x$.

\item[(\ref{axioma_A2})] $(x\too y)\land y=((x\to y)\land (\sim y\to \sim x))\land y=((x\to y)\land y)\land(\sim y\to\sim x)\underset{\rm\ref{lema_prop_Nelson} (\ref{tesis_1_2})}{=}y\land(\sim y\to\sim x)$ $\underset{\rm(\ref{proposicion2_3})}{=}y$.

\item[(\ref{axioma_A3})] 
$x\land(x\too y)=x\land(x\to y)\land(\sim y\to\sim x)\underset{\rm(N\ref{relacImplicaInfN})}{=}x\land(\sim x\lor y)\land(\sim y\to\sim x)=$\linebreak[4] $((x\land\sim x)\lor(x\land y))\land(\sim y\to\sim x)=((x\land\sim x)\land(\sim y\to\sim x))\lor((x\land y)\land(\sim y\to\sim x))$ $\underset{\rm(\ref{proposicion2_3})}{=}(x\land\sim x)\lor((x\land y)\land(\sim y\to\sim x))=x\land(\sim x\lor(y\land (\sim y\to\sim x)))\underset{\rm(\ref{proposicion2_3})}{=}x\land(\sim x\lor y)$ $\underset{\rm(N\ref{involucionN})}{=}x\land(\sim x\lor\sim\sim y)\underset{\rm(N\ref{deMenInfN})}{=}x\land\sim(x\land\sim y)$.

\item[(\ref{axioma_A4})] 
$x\too(y\land z)=(x\to(y\land z))\land(\sim(y\land z)\to\sim x)\underset{\rm\ref{lema_prop_Nelson} (\ref{distribInfADerN})}{=}(x\to y)\land(x\to z)\land(\sim(y\land z)\to\sim x)\underset{\rm(N\ref{deMenInfN})}{=}(x\to y)\land(x\to z)\land((\sim y\lor \sim z)\to\sim x)\underset{\rm\ref{lema_prop_Nelson} (\ref{tesis_1_9})}{=}(x\to y)\land(x\to z)\land(\sim y\to\sim x)\land(\sim z\to\sim x)=((x\to y)\land(\sim y\to\sim x))\land((x\to z)\land(\sim z\to\sim x))=(x\too y)\land(x\too z)$.
	
\item[(\ref{axioma_A5})] It is an immediate consequence of the definition of $\too$ and $\rm(N\ref{involucionN})$.
	
\item[(\ref{axioma_A6})]  $x\too(x\too(y\too(y\too z)))\underset{\rm\ref{lema_prop_Nelson} (\ref{tesis_1_83})}{=}x\to(y\to z)\underset{\rm(N\ref{relacImplicaInf2N})}{=}(x\land y)\to z\underset{\rm\ref{lema_prop_Nelson} (\ref{tesis_1_83})}{=}(x\land y)\too((x\land y)\too z)$.

\item[(\ref{axioma_A7})]  
Observe that $\sim(\sim x \land y)\too(x\too y)\geq x\too y$ follows from (\ref{axioma_A2}). Let us see the other inequality. From the definition of $\too$ , we can deduce the following: $u\too v\leq u\to v$.

Therefore, we have that 
$\sim(\sim x \land y)\too(x\too y)\leq
\sim(\sim x \land y)\to(x\too y)\underset{\text{def.}}{=}$ $
(x\lor\sim y)\to((x\to y)\land(\sim y\to\sim x))\underset{\rm\ref{lema_prop_Nelson} (\ref{distribInfADerN})}{=}
((x\lor\sim y)\to(x\to y))\land((x\lor\sim y)\to(\sim y\to\sim x))$.

Therefore
\begin{equation} \label{230919_01_I}
\sim(\sim x \land y)\too(x\too y)\leq ((x\lor\sim y)\to(x\to y))\land((x\lor\sim y)\to(\sim y\to\sim x)).
\end{equation}

Now, $(x\lor\sim y)\to(x\to y)\underset{\rm\ref{lema_prop_Nelson} (\ref{tesis_1_9})}{=}(x\to (x\to y))\land(\sim y\to(x\to y))$ \\ $\underset{\rm\ref{lema_prop_Nelson}(\ref{tesis_1_34})}{=}(x\to y)\land(\sim y\to(x\to y))$ $\leq$ $x\to y$. Hence,
\begin{equation} \label{230919_02_i}
(x\lor\sim y)\to(x\to y)\leq x\to y.
\end{equation}
On the other hand, $(x\lor\sim y)\to(\sim y\to\sim x)\underset{\rm\ref{lema_prop_Nelson}(\ref{tesis_1_9})}{=}(x\to(\sim y\to\sim x))\land(\sim y\to(\sim y\to\sim x))$ $\underset{\rm\ref{lema_prop_Nelson}(\ref{tesis_1_34})}{=}(x\to(\sim y\to \sim x))\land(\sim y\to\sim x)\underset{\rm\ref{lema_prop_Nelson} (\ref{tesis_1_2})}{=}\sim y\to\sim x$. Then,
\begin{equation} \label{230919_03_ii}
(x\lor\sim y)\to(\sim y\to\sim x)=\sim y\to\sim x.
\end{equation}

From (\ref{230919_01_I}), (\ref{230919_02_i}) and (\ref{230919_03_ii}) we conclude that\\ $\sim(\sim x \land y)\too(x\too y)\leq$ $(x\to y)\land(\sim y\to\sim x)=x\too y$.
\end{itemize}	
\end{proof}

\begin{Teo}\label{teodenelsonabrignole}
Let $\mathbf A=\langle A,\land,\lor,\to,\sim,1 \rangle$ be a Nelson algebra. Then $\langle A,\landb,\too,0 \rangle$ is a Brignole algebra.
\end{Teo}

\begin{proof}
From items (\ref{axioma_A1}), (\ref{axioma_A2}), (\ref{axioma_A3}), (\ref{axioma_A4}), (\ref{axioma_A5}), (\ref{axioma_A6}) and (\ref{axioma_A7}) of Lemma \ref{lema_prop_Nelson_Nuevas} it follows that $\mathbf A$ satisfies (B\ref{AxiomaBrignole01}), (B\ref{AxiomaBrignole02}), (B\ref{AxiomaBrignole03}), (B\ref{AxiomaBrignole04}), (B\ref{AxiomaBrignole05}), (B\ref{AxiomaBrignole06}) and (B\ref{AxiomaBrignole07}) respectively. The axioms (B\ref{AxiomaBrignole08}), (B\ref{AxiomaBrignole09}) and (B\ref{AxiomaBrignole10}) are verified considering Lemma \ref{lema_prop_Nelson_Nuevas} (\ref{axiomasReticuladoB}).
\end{proof}

\bigskip

Now, let us consider a Brignole algebra $\mathbf{A}=\langle A,\landb,\too,0 \rangle$. We define over $\mathbf{A}$ the following:

\begin{itemize}
	\item $x\land y:=x\landb y$,
	\item $x\lor y:=((x\too 0)\landb(y\too 0))\too 0$,
	\item $x\to y:=x\too(x\too y)$,
	\item $\sim x:=x\too 0$,
	\item $1:=0\too 0$.
\end{itemize}

Our goal now is to prove that $\langle A,\land,\lor,\to,\sim,1 \rangle$ is a Nelson algebra.

\begin{Lema} \label{lema_prop_Brignole}
	In a Brignole algebra $\mathbf A$ the following conditions are satisfied for all $x,y,z\in A$:
	\begin{enumerate}[\rm (a)]
		\item $x\land(x\lor y)=x$, $x\land(y\lor z)=(z\land x)\lor(y\land x)$ and $(x\land\sim x)\land(y\lor\sim y)=x\land\sim x$,\label{axiomasReticuladoyDMBrignole}
		\item $\sim 1=0$,\label{regla1}
		\item $1=1\too1$,\label{regla2}
		\item $1\too x=x$,\label{regla3}
		\item $\sim\sim x=x$,\label{AxN4}
		\item $\sim(x\lor y)=\sim x\land\sim y$,\label{regla4}
		\item $\sim(x\land y)=\sim x\lor \sim y$,\label{AxN5}
		\item $x\land (x\to y)=x\land(\sim x\lor y)$,\label{AxN8} 
		\item $x\lor 1=1$,\label{AxN1} 
		\item $x\too 1=1$,\label{regla5}
		\item $(x\land y)\to z=x\to(y\to z)$,\label{AxN9}
		\item $x\too x=1$, and in particular, $x\too x=y\too y$.\label{regla7}
		\item $x\to x=1$.\label{AxN7}
	\end{enumerate}
\end{Lema}

\begin{proof}
\begin{itemize}
\item[(\ref{axiomasReticuladoyDMBrignole})] 
It is an immediate consequence of (B\ref{AxiomaBrignole08}), (B\ref{AxiomaBrignole09}) and (B\ref{AxiomaBrignole10}). 	

From now on, we will use the fact that the reduct $\langle A,\land,\lor\rangle$ is a distributive lattice \cite{sholander51postulates}, with all its inherent properties.

\item[(\ref{regla1})]  $\sim 1=\sim(0\too 0)\underset{\text{def.}}{=}(0\too0)\too0\underset{\rm(B\ref{AxiomaBrignole01})}{=}0$.

\item[(\ref{regla2})] $1=0\too0\underset{\rm(B\ref{AxiomaBrignole05})}{=}\sim0\too\sim0=(0\too0)\too(0\too0)\underset{\text{def.}}{=}1\too1$.

\item[(\ref{regla3})]
 $1\too x\underset{\rm(\ref{regla2})}{=} (1\too 1)\too x\underset{\rm(B\ref{AxiomaBrignole01})}{=}x$.
 
\item[(\ref{AxN4})]
$x\underset{\rm(B\ref{AxiomaBrignole01})}{=}(y\too y)\too x\underset{\rm(B\ref{AxiomaBrignole05})}{=}
(x\too 0)\too((y\too y)\too 0)\underset{\rm(B\ref{AxiomaBrignole01})}{=}
(x\too 0)\too 0\underset{\rm\text{def.}~\sim}{=}\sim\sim x$.

\item[(\ref{regla4})]
 $\sim(x\lor y)=\sim(((x\too0)\land(y\too0))\too0)=\sim(\sim(\sim x\land\sim y))\underset{\rm(\ref{AxN4})}{=}\sim x\land\sim y$.
 
\item[(\ref{AxN5})] It follows from items (\ref{AxN4}) and (\ref{regla4}).
 
\item[(\ref{AxN8})]  $x\land (x\to y)\underset{\text{def.}}{=}x\land(x\too(x\too y))\underset{\rm(B\ref{AxiomaBrignole03})}{=}x\land\sim(x\land\sim(x\too y))\underset{\rm(\ref{AxN4})\ and\ (\ref{AxN5})}{=}x\land(\sim x\lor(x\too y))=(x\land\sim x)\lor(x\land(x\too y))\underset{\rm(B\ref{AxiomaBrignole03})}{=}(x\land\sim x)\lor(x\land\sim(x\land\sim y))=x\land(\sim x\lor\sim(x\land\sim y))\underset{\rm(\ref{AxN5})}{=}$ $x\land(\sim x\lor(\sim x\lor\sim\sim y))\underset{\rm(\ref{AxN4})}{=}x\land(\sim x\lor\sim x\lor y)=x\land(\sim x\lor y)$.

\item[(\ref{AxN1})] $x\lor 1=x\lor\sim 0\underset{\text{def.}}{=}\sim(\sim x\land\sim\sim0)\underset{\rm(\ref{AxN4})}{=}\sim(\sim x\land0)\underset{\text{def.}}{=}((x\too0)\land0)\too0\underset{\rm(B\ref{AxiomaBrignole02})}{=} 0\too0=1$.

By this result, we can conclude that 1 is the top element of $A$.

\item[(\ref{regla5})] $1\land(x\too 1)\underset{\rm(B\ref{AxiomaBrignole02})}{=}1$, then $1\leq x\too 1$. By (\ref{AxN1}), the equality follows.

\item[(\ref{AxN9})]
	By (B\ref{AxiomaBrignole06}), $(x\land y)\to z=(x\land y)\too((x\land y)\too z)=x\too(x\too(y\too(y\too z)))=x\to(y\to z)$.

\item[(\ref{regla7})]
$x\too x\underset{\rm(\ref{AxN4})}{=}
    ((x\too x)\too0)\too0)\underset{\rm(B\ref{AxiomaBrignole01})}{=}0\too0\underset{\text{def.}}{=}1$.
	
\item[(\ref{AxN7})] 
	$x\to x=x\too(x\too x)\underset{\rm(\ref{regla7})}{=}x\too1\underset{\rm(\ref{regla5})}{=}1$.
\end{itemize}
\end{proof}

\begin{Teo}\label{teodebrignoleanelson}
	Let $\langle A,\landb,\too,0\rangle$ be a Brignole algebra. Then $\mathbf A=\langle A,\land,\lor,\to,\sim,1\rangle$ is a Nelson algebra.
\end{Teo}

\begin{proof}
By Lemma \ref{lema_prop_Brignole} (\ref{axiomasReticuladoyDMBrignole}), $\mathbf A$ satisfies (N\ref{primercondicretN}), (N\ref{segundacondicretN}) and (N\ref{kleeneN}). The items (\ref{AxN4}), (\ref{AxN5}), (\ref{AxN7}), (\ref{AxN8}) and (\ref{AxN9}) from Lemma \ref{lema_prop_Brignole} prove the validity of (N\ref{involucionN}), (N\ref{deMenInfN}), (N\ref{ximpxigual1N}), (N\ref{relacImplicaInfN}) and (N\ref{relacImplicaInf2N}) respectively.
\end{proof}

\bigskip

\begin{Teo} The varieties of Nelson and Brignole algebras are term equivalent.
\end{Teo}

\begin{proof}
	If a Nelson algebra $\langle A,\land,\lor,\to,\sim,1 \rangle$ is obtained from a Brignole algebra $\langle A,\land,\too,0\rangle$ as in Theorem \ref{teodenelsonabrignole}, and if we define $x\Rightarrow y:=(x\to y)\land(\sim y\to \sim x)$, we obtain $x\Rightarrow y=x\too y$:

	$x\Rightarrow y=(x\to y)\land(\sim y\to \sim x)\underset{\text{def.}}{=}(x\too(x\too y))\land(\sim y\too(\sim y\too \sim x))\underset{\rm(B\ref{AxiomaBrignole05})}{=}\\
	(x\too(x\too y))\land(\sim y\too(x\too y))\underset{\rm(B\ref{AxiomaBrignole05})}{=}(\sim(x\too y)\too\sim x)\land(\sim(x\too y)\too y)\underset{\rm(B\ref{AxiomaBrignole04})}{=}\\
	\sim(x\too y)\too(\sim x\land y)\underset{\rm(B\ref{AxiomaBrignole05})}{=}\sim(\sim x\land y)\too(x\too y)\underset{\rm(B\ref{AxiomaBrignole07})}{=}x\too y$.

	If a Brignole algebra $\langle A,\landb,\too,0 \rangle$ is obtained from a Nelson algebra $\langle A,\land,\lor,\to,\sim,1 \rangle$ as in Theorem \ref{teodebrignoleanelson}, when we define $x\rightsquigarrow y:=x\too(x\too y)$, we obtain that $x\rightsquigarrow y=x\to y$. This is a consequence of Lemma \ref{lema_prop_Nelson} (\ref{tesis_1_83}).
\end{proof}

\section{Independence of Brignole axioms}

A natural question is which axioms of Definition \ref{defbrig} are independent. We have the following result:

\begin{Teo}
    In the variety $\Br$ the axioms $\rm(B1)$, $\rm(B3)$, $\rm(B5)$, $\rm(B7)$ and $\rm(B9)$ are independent.
\end{Teo}
\begin{proofish}
    The examples in this section have been found by the program Prover9 and Mace4 \cite{prover9}. For each example, we indicate the elements for which the equation fails, while the rest of them have been checked to hold.
\end{proofish}


\subsection{Independence of (B1)}
	\begin{center}
	\begin{minipage}{5cm}
	\begin{tabular}{c|ccc}
		$\too$ & 0 & $a$ & 1 \\ \hline
		      0 & 1 & 1 & 1 \\
		      $a$ & $a$ & $a$ & 1 \\
		      1 & 0 & $a$ & 1
	\end{tabular}
	\end{minipage}
	\begin{minipage}{5cm}
	\begin{tabular}{c|ccc}
		$\landb$ & 0 & $a$ &  1\\ \hline
		           0 & 0 & 0 & 0\\
		           $a$ & 0 & $a$ & $a$\\
		           1 & 0 & $a$ & 1
	\end{tabular}
	\end{minipage}
	\end{center}
\medskip
	
Axiom (B1) fails considering $x=a$ and $y=0$. Indeed: $(a\too a)\too 0=a\too0=a\neq0$.
	
\subsection{Independence of (B3)}
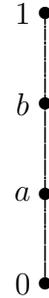
\begin{center}	
\begin{minipage}{5cm}
	\begin{tabular}{r|rrrr}
		$\too$ & 0 & $a$ & $b$ & 1\\
		\hline
		0 & 1 & 1 & 1 & 1 \\
		$a$ & $b$ & 1 & 1 & 1 \\
		$b$ & $a$ & 1 & 1 & 1 \\
		1 & 0 & $a$ & $b$ & 1
	\end{tabular}
\end{minipage}
\begin{minipage}{5cm}
	\begin{tabular}{r|rrrr}
		$\landb$ & 0 & $a$ & $b$ & 1\\
		\hline
		0 & 0 & 0 & 0 & 0 \\
		$a$ & 0 & $a$ & $a$ & $a$ \\
		$b$ & 0 & $a$ & $b$ & $b$ \\
		1 & 0 & $a$ & $b$ & 1
	\end{tabular} 
\end{minipage}
\begin{minipage}{3cm}
	\beginpicture
	\setcoordinatesystem units <3mm,3mm>
	\setplotarea x from 0 to 16, y from 0 to 8
	\put {$ \bullet$} [c] at  8 0
	\put {$ 0$} [c] at  7 0
	\put {$ \bullet$} [c] at 8 4
	\put {$a$} [c] at  7 4
	\put {$ \bullet$} [c] at  8 8
	\put {$b$} [c] at  7 8
	\put {$ \bullet$} [c] at  8 12
	\put {$1$} [c] at  7 12
	\setlinear \plot 8 0 8 12 /
	\endpicture
\end{minipage}
\end{center}
	
Axiom (B3) fails considering $x=a$ and $y=b$.
	
Indeed: $a\landb((a\landb(b\too0))\too0)=a\landb((a\landb0)\too0)=a\landb(0\too0)=a\landb0=0$, and $a\landb(a\too b)=a\landb1=a$.
	
\subsection{Independence of (B5)}
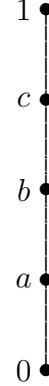
\begin{center}	
	\begin{minipage}{5cm}
	\begin{tabular}{r|rrrrr}
		$\too$ & 0 & $a$ & $b$ & $c$ & 1 \\ \hline
		      0 & 1 & 1 & 1 & 1 & 1 \\
		      $a$ & $c$ & 1 & 1 & 1 & 1 \\
		      $b$ & $b$ & 1 & 1 & 1 & 1 \\
		      $c$ & $a$ & $a$ & $b$ & 1 & 1 \\
		      1 & 0 & $a$ & $b$ & $c$ & 1
	\end{tabular} 
\end{minipage}
\hspace{.5cm}
	\begin{minipage}{5cm}
	\begin{tabular}{r|rrrrr}
		$\landb$ & 0 & $a$ & $b$ & $c$ & 1 \\ \hline
		       0 & 0 & 0 & 0 & 0 & 0 \\
		       $a$ & 0 & $a$ & $a$ & $a$ & $a$ \\
		       $b$ & 0 & $a$ & $b$ & $b$ & $b$ \\
		       $c$ & 0 & $a$ & $b$ & $c$ & $c$ \\
		       1 & 0 & $a$ & $b$ & $c$ & 1
	\end{tabular}
\end{minipage}
\begin{minipage}{5cm}
\beginpicture
\setcoordinatesystem units <3mm,3mm>
\setplotarea x from 0 to 16, y from 0 to 8
\put {$ \bullet$} [c] at  8 0
\put {$ 0$} [c] at  7 0
\put {$ \bullet$} [c] at 8 4
\put {$a$} [c] at  7 4
\put {$ \bullet$} [c] at  8 8
\put {$b$} [c] at  7 8
\put {$ \bullet$} [c] at  8 12
\put {$c$} [c] at  7 12
\put {$ \bullet$} [c] at  8 16
\put {$1$} [c] at  7 16
\setlinear \plot 8 0 8 16 /
\endpicture
\end{minipage}
\end{center}
	
Axiom (B5) fails considering $x=c$ and $y=b$.
	
Indeed: $c\too b=b$, and $(b\too 0)\too(c\too0)=b\too a=1$.

\subsection{Independence of (B6)}
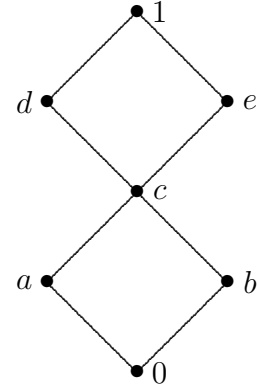
\begin{center}	
	\begin{minipage}{4.8cm}
	\begin{tabular}{r|rrrrrrr}
		$\too$ & 0 & $a$ & $b$ & $c$ & $d$ & $e$ & 1\\ \hline
		0 & 1 & 1 & 1 & 1 & 1 & 1 & 1 \\
		$a$ & $e$ & 1 & $e$ & 1 & 1 & 1 & 1 \\
		$b$ & $d$ & $d$ & 1 & 1 & 1 & 1 & 1 \\
		$c$ & $c$ & $d$ & $e$ & 1 & 1 & 1 & 1 \\
		$d$ & $b$ & $c$ & $b$ & $e$ & 1 & $e$ & 1 \\
		$e$ & $a$ & $a$ & $c$ & $d$ & $d$ & 1 & 1 \\
		1 & 0 & $a$ & $b$ & $c$ & $d$ & $e$ & 1
	\end{tabular}
\end{minipage}
 \hspace{.5cm}
 \begin{minipage}{4.8cm}
    \begin{tabular}{r|rrrrrrr}
	$\landb$ & 0 & $a$ & $b$ & $c$ & $d$ & $e$ & 1\\ \hline
	0 & 0 & 0 & 0 & 0 & 0 & 0 & 0 \\
	$a$ & 0 & $a$ & 0 & $a$ & $a$ & $a$ & $a$ \\
	$b$ & 0 & 0 & $b$ & $b$ & $b$ & $b$ & $b$ \\
	$c$ & 0 & $a$ & $b$ & $c$ & $c$ & $c$ & $c$ \\
	$d$ & 0 & $a$ & $b$ & $c$ & $d$ & $c$ & $d$ \\
	$e$ & 0 & $a$ & $b$ & $c$ & $c$ & $e$ & $e$ \\
	1 & 0 & $a$ & $b$ & $c$ & $d$ & $e$ & 1
	\end{tabular}
\end{minipage}
	\begin{minipage}{5cm}
	\beginpicture
	\setcoordinatesystem units <3mm,3mm>
	\setplotarea x from 0 to 16, y from 0 to 8
	\put {$ \bullet$} [c] at  8 0
	\put {$ 0$} [c] at  9 0
	\put {$ \bullet$} [c] at 4 4
	\put {$a$} [c] at  3 4
	\put {$ \bullet$} [c] at  12 4
	\put {$b$} [c] at  13 4
	\put {$ \bullet$} [c] at 8  8 
	\put {$c$} [c] at  9 8
	\put {$ \bullet$} [c] at  4 12
	\put {$d$} [c] at  3 12
	\put {$ \bullet$} [c] at  12 12
	\put {$e$} [c] at 13 12
	\put {$ \bullet$} [c] at  8 16
	\put {$1$} [c] at  9 16
	
	\setlinear \plot 8 0 4 4 12 12 8 16 4 12 12 4 8 0 /
	\endpicture
\end{minipage}
\end{center}

\medskip

Axiom (B6) fails considering $x=e$, $y=d$ and $z=0$.
	
Indeed: $e\too(e\too(d\too(d\too0)))=$ $e\too(e\too(d\too b))=e\too(e\too b)=e\too c=d$, and $(e\landb d)\too((e\landb d)\too0)=c\too(c\too0)=c\too c=1$.

\subsection{Independence of (B9)}
\begin{center}
 \begin{minipage}{5cm}
	\begin{tabular}{r|rrrr}
		$\too$ & 0 & $a$ & $b$ & 1 \\ \hline
		      0 & 1 & 1 & 1 & 1 \\
		      $a$ & $a$ & 1 & 1 & 1 \\
		      $b$ & $b$ & 1 & 1 & 1 \\
		      1 & 0 & $a$ & $b$ & 1
	\end{tabular}
\end{minipage}
 \hspace{.5cm}
\begin{minipage}{5cm}
	\begin{tabular}{r|rrrr}
		$\landb$ & 0 & $a$ & $b$ & 1 \\ \hline
		        0 & 0 & 0 & 0 & 0 \\
		        $a$ & 0 & $a$ & $a$ & $a$ \\
		        $b$ & 0 & $b$ & $b$ & $b$ \\
		        1 & 0 & $a$ & $b$ & 1
	\end{tabular}
\end{minipage}
\end{center}
	
Axiom (B9) fails considering $x=b$, $y=0$ and $z=a$.
	
Indeed: $b\landb(0\lorb a)=b\landb(((0\too0)\landb(a\too0))\too0)=b\landb((1\landb a)\too0)=b\landb(a\too0)=b\landb a=b$, and $(a\landb b)\lorb (0\landb b)=a\lorb 0=((a\too0)\landb(0\too0))\too0=(a\landb1)\too0=a\too0=a$.\hfill $\blacksquare$

\section{Dependent axioms}
In this section we will prove that axioms (B2) and (B8) can be derived from the other axioms for Brignole algebras.

\begin{Lema} \label{propextrabrignole}
	Let $\bf A$ be an algebra  $\langle A,\landb,\too,0\rangle$ satisfying the axioms {\rm (B1), (B3)} to {\rm(B7), (B9)}, and {\rm (B10)}. The following properties are satisfied for all $x,y,z\in A$:
    \begin{enumerate}[\rm (a)]
        \item $\sim 1=0$,\label{28052021_1}
    	\item $\sim\sim x=x$,\label{AxN4_v}
    	\item $x\too x=1$, and in particular, $x\too x=y\too y$,\label{28052021}
    	\item $x=(x\too0)\too(x\landb0)$,\label{propiedadextraLuiz}
 		\item $0\landb0=0$,\label{infceroceroescero}
 		\item $x\landb y=y\landb x$,\label{conmutinfimo}
    	\item $x\landb x=x$,\label{idempotenciainfimo}
    	\item $0\landb 1=0$,\label{infimoceroconuno}
    	\item $x\landb 0=0$,\label{infimoconcero}
        \item $x\lorb x=x$,\label{idempsupremo}
        \item $x\lorb y=y\lorb x$,\label{conmutsupremo}
        \item $\sim(x\lorb y)=\sim x\landb\sim y$, and $\sim(x\landb y)=\sim x\lorb\sim y$,\label{demorganbrignole}
        \item $x\landb(y\lorb z)=(x\landb y)\lorb(x\landb z)$,\label{distribBrignole}
        \item $x\lorb 1=1$.\label{supremoconuno}
    \end{enumerate}
\end{Lema}

\begin{proof}   	\setcounter{equation}{0}
	Items (\ref{28052021_1}), (\ref{AxN4_v}) and (\ref{28052021}) were proved in Lemma \ref{lema_prop_Brignole} (\ref{regla1}), (\ref{AxN4}) and (\ref{regla7}) respectively, without making use of (B2) nor (B8).
	\begin{itemize}
	\item[(\ref{propiedadextraLuiz})]
	    $x\underset{\rm(\ref{AxN4_v})}{=}
    	(x\too0)\too0=
    	\sim x\too0\underset{\rm(B\ref{AxiomaBrignole07})}{=}\sim(\sim\sim x\landb0)\too(\sim x\too0)\underset{\rm(\ref{AxN4_v})\text{ and def.~}\sim}{=}\sim(x\landb0)\too\sim\sim x\underset{\rm(B\ref{AxiomaBrignole05})}{=}$ $\sim x\too(x\landb0)=(x\too0)\too(x\landb0)$.
    \item[(\ref{infceroceroescero})] 
    	Taking $x$ to be $0$ in (\ref{propiedadextraLuiz}), we have that $0=(0\too0)\too(0\landb0)\underset{\rm(B\ref{AxiomaBrignole01})}{=}0\landb0$.
	\item[(\ref{conmutinfimo})]
$x\landb y\underset{\rm(\ref{AxN4_v})}{=}((x\landb y)\too0)\too0\underset{\rm(\ref{infceroceroescero})}{=}((x\landb y)\too(0\landb0))\too0\underset{\rm(B\ref{AxiomaBrignole04})}{=}$\\
	    $(((x\landb y)\too0)\landb((x\landb y)\too0))\too0\underset{\rm def.~}{=}(x\landb y)\lorb(x\landb y)\underset{\rm(B\ref{AxiomaBrignole09})}{=}
	    y\landb(x\lorb x)\underset{\rm def.~}{=}$\\$y\landb(((x\too0)\landb(x\too0))\too0)\underset{\rm(B\ref{AxiomaBrignole04})}{=}$
	    $y\landb((x\too(0\landb0))\too0)\underset{\rm(\ref{infceroceroescero})}{=}y\landb((x\too0)\too0)\underset{\rm(\ref{AxN4_v})}{=}y\landb x.$


	 \item[(\ref{idempotenciainfimo})]
	    $x\landb x\underset{\rm(\ref{AxN4_v})}{=}\sim\sim x \landb \sim \sim x\underset{\rm def.~}{=}((x\too0)\too0)\landb((x\too0)\too0)\underset{\rm (B\ref{AxiomaBrignole04})}{=}(x\too0)\too(0\landb 0)\underset{\rm (\ref{infceroceroescero})}{=}$ $(x\too 0)\too0\underset{\rm def.~}{=}\sim\sim x\underset{\rm(\ref{AxN4_v})}{=}x$.
	 
	 
    	    
    \item[(\ref{infimoceroconuno})] 
    
    We notice that $x\underset{\rm (\ref{propiedadextraLuiz})}{=}(x\too0)\too(x\landb0)\underset{\rm (\ref{conmutinfimo}), def.}{=}\sim x\too(0\landb x)$.
        
    Replacing $x$ with $\sim y$, we obtain the equivalent
        \begin{equation}\label{30052021_a}
            \sim y=y\too(0\landb \sim y).
        \end{equation}
        
        Using (B\ref{AxiomaBrignole03}) and the definition of $\sim$, we have that $0\landb(0\too x)=0\landb((0\landb(x\too0))\too0)$. Then
        \begin{equation}\label{30052021_c}
            (0\landb(x\too 0))\too(0\landb(0\too x))=(0\landb(x\too0))\too(0\landb(0\landb(x\too0))\too0).
        \end{equation}
        The right side of (\ref{30052021_c}) is of the same form as the right side in (\ref{30052021_a}) taking $y$ to be $0\landb(x\too0)$, so we can rewrite (\ref{30052021_c}) as
        \begin{equation}\label{30052021_d}
            (0\landb(x\too 0))\too(0\landb(0\too x))=(0\landb(x\too0))\too0.
        \end{equation}
	    Replacing $x$ by 0 in (\ref{30052021_d}), we obtain
	    $$(0\landb(0\too0))\too0=(0\landb(0\too0))\too(0\landb(0\too 0))\underset{\rm(\ref{28052021})}{=}y\too y\underset{\rm(\ref{28052021})}{=}0\too 0,$$
	    that is
	    \begin{equation}\label{30052021_e}
	        (0\landb(0\too0))\too0=0\too 0.
	    \end{equation}
	    By (\ref{AxN4_v}), (\ref{30052021_e}) is equivalent to
	    $$0\landb(0\too 0)=0.$$
	    Therefore, $0\landb 1=0$.
	    \item[(\ref{infimoconcero})]
	    See the Appendix.
	
    \item [(\ref{idempsupremo})] $x\lorb x\underset{\rm def.~}{=}((x\too 0)\landb(x\too0))\too0\underset{\rm(\ref{idempotenciainfimo})}{=}(x\too 0)\too0\underset{\rm(\ref{AxN4_v})}{=}x$.
    \item [(\ref{conmutsupremo})] $x\lorb y\underset{\rm def.~}{=}((x\too 0)\landb(y\too0))\too0\underset{\rm(\ref{conmutinfimo})}{=}((y\too 0)\landb(x\too0))\too0\underset{\rm def.~}{=}y\lorb x$.
    \item [(\ref{demorganbrignole})] This was proved in Lemma \ref{lema_prop_Brignole} (\ref{regla4}) and (\ref{AxN5}), without use of (B2) nor (B8).
    \item [(\ref{distribBrignole})] $x\landb(y\lorb z)\underset{\rm (B\ref{AxiomaBrignole09})}{=}(z\landb x)\lorb(y\landb x)\underset{\rm(\ref{conmutinfimo})}{=}(x\landb z)\lorb(x\landb y)\underset{\rm(\ref{conmutsupremo})}{=}(x\landb y)\lorb(x\landb z)$.
    \item [(\ref{supremoconuno})] $x\lorb 1\underset{\rm(\ref{AxN4_v})}{=}\sim\sim x\lorb\sim\sim1\underset{\rm(\ref{demorganbrignole})}{=}\sim(\sim x\landb\sim1)\underset{\rm(\ref{28052021_1})}{=}\sim(\sim x\landb 0)\underset{\rm(\ref{infimoconcero})}{=}\sim 0\underset{\rm def.~}{=}1$.
    \end{itemize}
\end{proof}

We are now in condition of showing the following:

\begin{Teo}
	An algebra $\langle A,\landb,\too,0\rangle$ of type $(2,2,0)$ is a Brignole algebra if and only if it satisfies the following equations for every $x,y,z\in A$:
	\begin{itemize}\label{teoIndep}
		\item[\rm (B1)] $(x\too x)\too y=y$,
		\item[\rm(B3)] $x\landb\simb(x\landb\simb y)=x\landb(x\too y)$,
		\item[\rm(B4)] $x\too(y\landb z)=$ $(x\too y)\landb(x\too z)$,
		\item[\rm(B5)] $x\too y=\simb y\too\simb x$,
		\item[\rm(B6)] $x\too(x\too(y\too(y\too z)))=(x\landb y)\too((x\landb y)\too z)$,
		\item[\rm(B7)] $\simb(\simb x\landb y)\too(x\too y)=x\too y$,
		\item[\rm(B9)] $x\landb(y\lorb z)=(z\landb x)\lorb(y\landb x)$,
		\item[\rm(B10)] $(x\landb\simb x)\landb (y\lorb\simb y)=x\landb\simb x$,
	\end{itemize}
	where   $\simb x:=x\too 0$ and  $x\lorb y:=((x\too 0)\landb(y\too 0))\too 0$.
\end{Teo}
\begin{proof}
	\setcounter{equation}{0}
    One implication is immediate. For the other one, let us prove (B2) first. The letters in the proof reference items from Lemma \ref{propextrabrignole}.
    
    $(x\too y)\landb y\underset{\rm(B5)}{=}((y\too0)\too(x\too0))\landb y\underset{\rm (\ref{AxN4_v})}{=}((y\too0)\too(x\too0))\landb(y\too0)\too0$ $\underset{\rm(B4)}{=}(y\too0)\too((x\too0)\landb 0)\underset{\rm  (\ref{infimoconcero})}{=}(y\too0)\too0\underset{\rm (\ref{AxN4_v})}{=}y$.
    
    Hence, we have (B2).
    
    Using (\ref{demorganbrignole}), (\ref{distribBrignole}), (\ref{conmutsupremo}), and (\ref{conmutinfimo}) of Lemma \ref{propextrabrignole}, we have that:
    $$\sim(x\landb\sim y)\landb\sim(z\landb\sim y)=\sim(\sim y\landb(z\lorb x)).$$
    
    Taking $z=1$ we obtain
    
    \begin{equation}\label{30052021}
        \sim(x\landb\sim y)\landb\sim(1\landb\sim y)=\sim(\sim y\landb(1\lorb x)).
    \end{equation}
    
    Notice that
    \[z\underset{\rm (B\ref{AxiomaBrignole02})}{=}z\landb(z\too z)\underset{\rm (\ref{28052021})}{=}z\landb 1\underset{\rm (\ref{conmutinfimo})}{=}1\landb z,\]
    
    that is
    
    \begin{equation}\label{30052021_1}
        z=1\landb z.
    \end{equation}
    
    Taking $z=\sim y$ in (\ref{30052021_1}), and replacing that in (\ref{30052021}), we obtain
    
    $$\sim(x\landb\sim y)\landb\sim\sim y=\sim(\sim y\landb (1\lorb x))\underset{\rm (\ref{conmutsupremo}),(\ref{supremoconuno})}{=}\sim(\sim y\landb 1)\underset{\rm (\ref{30052021_1}),(\ref{conmutinfimo})}{=}\sim\sim y\underset{\rm (\ref{AxN4_v})}{=}y,$$
    
    by Lemma \ref{propextrabrignole} (\ref{AxN4_v}) we have that
    $$\sim(x\landb\sim y)\landb y=y,$$
    
    and by Lemma \ref{propextrabrignole}  (\ref{distribBrignole}), (\ref{AxN4_v}) and (\ref{conmutinfimo}) it follows that
    
    \begin{equation}\label{30052021_2}
        y= y\landb(\sim x\lorb y).
    \end{equation}
    
    
    If we change $y$ and $x$ in (\ref{30052021_2}) by $x$ and $\sim y$ respectively (and use Lemma \ref{propextrabrignole} $\rm(\ref{conmutsupremo})$), we obtain $\rm(B8)$.
\end{proof}

\nocite{monteiro07works}
\bibliographystyle{alpha}

\section*{Appendix}
The following proof has been adapted from the output produced by the program {\sc Prover9}, \cite{prover9}. For simplicity we are going to replace $\landb$ with $\land$.


\begin{enumerate}[1.]

\item $x \lor_\text{\tiny B} y = ((x \too 0) \wedge (y \too 0)) \too 0$ \hfill definition of $\lorb$  \label{paso2}

\item $\sim x = x \too 0$ \hfill definition of $\sim$ \label{paso3}

\item $(x \too x) \too y = y$ \hfill $\rm(B1)$  \label{paso4}

\item $x \wedge \sim (x \wedge \sim y) = x \wedge (x \too y)$  \hfill $\rm(B3)$ \label{paso5}

\item $x \wedge ((x \wedge (y \too 0)) \too 0) = x \wedge (x \too y)$ \hfill by  (\ref{paso3}) and (\ref{paso5})  \label{paso7}

\item $x \too (y \wedge z) = (x \too y) \wedge (x \too z)$  \hfill $\rm(B4)$ \label{paso8}

\item $(x \too y) \wedge (x \too z) = x \too (y \wedge z)$ \hfill by (\ref{paso8})  \label{paso9}

\item $x \too y = \sim y \too \sim x$ \hfill $\rm(B5)$  \label{paso10}

\item $x \too y = (y \too 0) \too (x \too 0)$ \hfill by  (\ref{paso3}) and (\ref{paso10})  \label{paso13}

\item $x \too (x \too (y \too (y \too z))) = (x \wedge y) \too ((x \wedge y) \too z)$ \hfill $\rm(B6)$ \label{paso14}

\item $\sim (\sim x \wedge y) \too (x \too y) = x \too y$ \hfill $\rm(B7)$  \label{paso16}

\item $\sim ((x \too 0) \wedge y) \too (x \too y) = x \too y$ \hfill by  (\ref{paso3}) and (\ref{paso16})  \label{paso17}

\item $(((x \too 0) \wedge y) \too 0) \too (x \too y) = x \too y$ \hfill by (\ref{paso3}) and (\ref{paso17})  \label{paso18}

\item $x \wedge (y \lor_\text{\tiny B} z) = (z \wedge x) \lor_\text{\tiny B} (y \wedge x)$ \hfill $\rm(B9)$  \label{paso19}

\item $x \wedge (((y \too 0) \wedge (z \too 0)) \too 0) = (((z \wedge x) \too 0) \wedge ((y \wedge x) \too 0)) \too 0$ \hfill by (\ref{paso2}) and (\ref{paso19})  \label{paso21}

\item $x \too x = y \too y$ \hfill by Lemma \ref{lema_prop_Brignole} (\ref{regla7}) \label{paso22}

\item $x \wedge y = y \wedge x$ \hfill by Lemma \ref{propextrabrignole}  \rm(\ref{conmutinfimo}) \label{paso23}

\item $x \wedge x = x$ \hfill by Lemma \ref{propextrabrignole}  \rm(\ref{idempotenciainfimo}) \label{paso24}

\item $0 = 0 \wedge (0 \too 0)$ \hfill by Lemma \ref{propextrabrignole}  (\ref{infimoceroconuno}) \label{paso182}

\item $x \too 0 = x \too (0 \wedge (x \too 0))$ \hfill by Lemma (\ref{propextrabrignole}), (\ref{infimoceroconuno}), item (1)  \label{paso51}

\item $x \wedge ((x \wedge 0) \too 0) = x \wedge (x \too (y \too y))$ \hfill by  (\ref{paso4}) and (\ref{paso7})  \label{paso27}

\item $(x \too y) \wedge ((x \too (y \wedge 0)) \too 0) = (x \too y) \wedge ((x \too y) \too x)$ \hfill by  (\ref{paso9}) and (\ref{paso7})  \label{paso28}

\item $(x \too 0) \too ((y \too y) \too 0) = x$ \hfill by  (\ref{paso13}) and (\ref{paso4})  \label{paso29}

\item $(x \too 0) \too 0 = x$ \hfill by  (\ref{paso4}) and (\ref{paso29})  \label{paso30}

\item $0 \too (x \too 0) = x \too (y \too y)$ \hfill by  (\ref{paso4}) and (\ref{paso13})  \label{paso31}

\item $x \wedge (x \too y) = x \wedge ((y \too 0) \too (x \too 0))$ \hfill by  (\ref{paso13}) \label{paso33}

\item $(x \too y) \wedge ((y \too 0) \too z) = (y \too 0) \too ((x \too 0) \wedge z)$ \hfill by  (\ref{paso13}) and (\ref{paso9})  \label{paso35}

\item $((x \too 0) \too y) \wedge (z \too x) = (x \too 0) \too (y \wedge (z \too 0))$ \hfill by  (\ref{paso13}) and (\ref{paso9})  \label{paso36}

\item $(x \too (y \wedge z)) \too (((x \too y) \wedge (x \too z)) \too u) = (x \too y) \too ((x \too y) \too ((x \too z) \too ((x \too z) \too u)))$ \hfill by  (\ref{paso9}) and (\ref{paso14})  \label{paso37}

\item $(x \too (y \wedge z)) \too ((x \too (y \wedge z)) \too u) = (x \too y) \too ((x \too y) \too ((x \too z) \too ((x \too z) \too u)))$ \\ \null  \hfill by  (\ref{paso9}) and (\ref{paso37})  \label{paso38}

\item $(x \wedge y) \too (x \too (x \too (y \too (y \too z)))) = x \too (x \too (y \too (y \too ((x \wedge y) \too z))))$ \\ \null  \hfill by  (\ref{paso14})   \label{paso39}

\item $((x \wedge y) \wedge (x \wedge y)) \too (((x \wedge y) \wedge (x \wedge y)) \too z) = (x \wedge y) \too (x \too (x \too (y \too (y \too ((x \wedge y) \too z)))))$ \\ \null\hfill by  (\ref{paso14})  \label{paso40}

\item $(x \wedge y) \too (((x \wedge y) \wedge (x \wedge y)) \too z) = (x \wedge y) \too (x \too (x \too (y \too (y \too ((x \wedge y) \too z)))))$\\\null \hfill by  (\ref{paso24}) and (\ref{paso40})  \label{paso41}

\item $(x \wedge y) \too ((x \wedge y) \too z) = (x \wedge y) \too (x \too (x \too (y \too (y \too ((x \wedge y) \too z)))))$ \hfill by  (\ref{paso24}) and (\ref{paso41})  \label{paso42}

\item $x \too (x \too (y \too (y \too z))) = (x \wedge y) \too (x \too (x \too (y \too (y \too ((x \wedge y) \too z)))))$ \\\null\hfill by  (\ref{paso14}) and (\ref{paso42})  \label{paso43}

\item $x \too (x \too (y \too (y \too z))) = x \too (x \too (y \too (y \too ((x \wedge y) \too ((x \wedge y) \too z)))))$ \\\null\hfill by  (\ref{paso39}) and (\ref{paso43})  \label{paso44}

\item $x \too (x \too (y \too (y \too z))) = x \too (x \too (y \too (y \too (x \too (x \too (y \too (y \too z)))))))$ \\\null \hfill by  (\ref{paso14}) and (\ref{paso44})  \label{paso45}

\item $((x \too y) \too 0) \too ((((x \too 0) \wedge y) \too 0) \too 0) = x \too y$ \hfill by  (\ref{paso18}) and (\ref{paso13})  \label{paso53}

\item $((x \too y) \too 0) \too ((x \too 0) \wedge y) = x \too y$ \hfill by  (\ref{paso30}) and (\ref{paso53})  \label{paso54}

\item $x \lor_\text{\tiny B} ((y \too 0) \wedge (z \too 0)) = ((((z \wedge (x \too 0)) \too 0) \wedge ((y \wedge (x \too 0)) \too 0)) \too 0) \too 0$ \\\null\hfill by  (\ref{paso21}) and (\ref{paso2})  \label{paso55}

\item $((x \too 0) \wedge (((y \too 0) \wedge (z \too 0)) \too 0)) \too 0 = ((((z \wedge (x \too 0)) \too 0) \wedge ((y \wedge (x \too 0)) \too 0)) \too 0) \too 0$ \hfill by  (\ref{paso2}) and (\ref{paso55})  \label{paso56}

\item $((x \too 0) \wedge (((y \too 0) \wedge (z \too 0)) \too 0)) \too 0 = ((z \wedge (x \too 0)) \too 0) \wedge ((y \wedge (x \too 0)) \too 0)$ \\\null\hfill by  (\ref{paso30}) and (\ref{paso56})  \label{paso57}

\item $((x \wedge (y \too 0)) \too 0) \wedge (((z \too 0) \wedge (x \too 0)) \too 0) = (((x \wedge (x \too y)) \too 0) \wedge ((z \wedge ((x \wedge (y \too 0)) \too 0)) \too 0)) \too 0$ \hfill by  (\ref{paso7}) and (\ref{paso21})  \label{paso60}

\item $(x \wedge (((y \too 0) \wedge (z \too 0)) \too 0)) \wedge ((((z \wedge x) \too 0) \wedge ((y \wedge x) \too 0)) \too u) = (((z \wedge x) \too 0) \wedge ((y \wedge x) \too 0)) \too (0 \wedge u)$ \hfill by  (\ref{paso21}) and (\ref{paso9})  \label{paso61}

\item $(x \too x) \wedge (y \too z) = y \too (y \wedge z)$ \hfill by  (\ref{paso22}) and (\ref{paso9})  \label{paso65}

\item $(x \too y) \wedge (z \too z) = x \too (y \wedge x)$ \hfill by  (\ref{paso22}) and (\ref{paso9})  \label{paso66}

\item $(x \wedge y) \too (z \too z) = x \too (x \too (y \too (y \too (x \wedge y))))$ \hfill by  (\ref{paso22}) and (\ref{paso14})  \label{paso67}

\item $0 \too ((x \wedge y) \too 0) = x \too (x \too (y \too (y \too (x \wedge y))))$ \hfill by  (\ref{paso31}) and (\ref{paso67})  \label{paso68}

\item $(((x \too 0) \wedge x) \too 0) \too (y \too y) = x \too x$ \hfill by  (\ref{paso22}) and (\ref{paso18})  \label{paso69}

\item $((x \wedge (x \too 0)) \too 0) \too (y \too y) = x \too x$ \hfill by  (\ref{paso23}) and (\ref{paso69})  \label{paso70}

\item $0 \too (((x \wedge (x \too 0)) \too 0) \too 0) = x \too x$ \hfill by  (\ref{paso31}) and (\ref{paso70})  \label{paso71}

\item $0 \too (x \wedge (x \too 0)) = x \too x$ \hfill by  (\ref{paso30}) and (\ref{paso71})  \label{paso72}

\item $x \wedge (((y \too 0) \wedge x) \too 0) = x \wedge (x \too y)$ \hfill by  (\ref{paso23}) and (\ref{paso7})  \label{paso73}

\item $(x \too 0) \wedge ((x \too 0) \too 0) = (x \too 0) \wedge ((x \too 0) \too x)$ \hfill by  (\ref{paso24}) and (\ref{paso7})  \label{paso74}

\item $(x \too 0) \wedge x = (x \too 0) \wedge ((x \too 0) \too x)$ \hfill by  (\ref{paso30}) and (\ref{paso74})  \label{paso75}

\item $x \wedge (x \too 0) = (x \too 0) \wedge ((x \too 0) \too x)$ \hfill by  (\ref{paso23}) and (\ref{paso75})  \label{paso76}

\item $x \too ((x \wedge x) \too y) = x \too (x \too (x \too (x \too y)))$ \hfill by  (\ref{paso24}) and (\ref{paso14})  \label{paso78}

\item $x \too (x \too y) = x \too (x \too (x \too (x \too y)))$ \hfill by  (\ref{paso24}) and (\ref{paso78})  \label{paso79}

\item $((x \too 0) \too 0) \too (x \too (x \too 0)) = x \too (x \too 0)$ \hfill by  (\ref{paso24}) and (\ref{paso18})  \label{paso81}

\item $x \too (x \too (x \too 0)) = x \too (x \too 0)$ \hfill by  (\ref{paso30}) and (\ref{paso81})  \label{paso82}

\item $0 \wedge (0 \too 0) = 0 \wedge (0 \too (x \too x))$ \hfill by  (\ref{paso24}) and (\ref{paso27})  \label{paso83}

\item $x \wedge ((x \wedge y) \too 0) = x \wedge (x \too (y \too 0))$ \hfill by  (\ref{paso30}) and (\ref{paso7})  \label{paso85}

\item $x \wedge ((x \too 0) \too y) = (x \too 0) \too (0 \wedge y)$ \hfill by  (\ref{paso30}) and (\ref{paso9})  \label{paso86}

\item $((x \too 0) \too y) \wedge x = (x \too 0) \too (y \wedge 0)$ \hfill by  (\ref{paso30}) and (\ref{paso9})  \label{paso87}

\item $x \wedge ((x \too 0) \too y) = (x \too 0) \too (y \wedge 0)$ \hfill by  (\ref{paso23}) and (\ref{paso87})  \label{paso88}

\item $x \too (y \too 0) = y \too (x \too 0)$ \hfill by  (\ref{paso30}) and (\ref{paso13})  \label{paso89}

\item $(x \too 0) \too y = (y \too 0) \too x$ \hfill by  (\ref{paso30}) and (\ref{paso13})  \label{paso90}

\item $((((x \too 0) \too 0) \wedge 0) \too 0) \too x = (x \too 0) \too 0$ \hfill by  (\ref{paso30}) and (\ref{paso18})  \label{paso91}

\item $((x \wedge 0) \too 0) \too x = (x \too 0) \too 0$ \hfill by  (\ref{paso30}) and (\ref{paso91})  \label{paso92}

\item $(x \too 0) \too (x \wedge 0) = (x \too 0) \too 0$ \hfill by  (\ref{paso90}) and (\ref{paso92})  \label{paso93}

\item $(x \too 0) \too (x \wedge 0) = x$ \hfill by  (\ref{paso30}) and (\ref{paso93})  \label{paso94}

\item $(x \too 0) \too (0 \wedge x) = x$ \hfill by  (\ref{paso94}) and (\ref{paso23})  \label{paso95}

\item $((x \too 0) \too x) \wedge ((x \too 0) \too 0) = x$ \hfill by  (\ref{paso94}) and (\ref{paso9})  \label{paso96}

\item $((x \too 0) \too x) \wedge x = x$ \hfill by  (\ref{paso30}) and (\ref{paso96})  \label{paso97}

\item $x \wedge ((x \too 0) \too x) = x$ \hfill by  (\ref{paso23}) and (\ref{paso97})  \label{paso98}

\item $((x \too 0) \too 0) \wedge (((x \too 0) \too (0 \wedge 0)) \too 0) = ((x \too 0) \too 0) \wedge (x \too (x \too 0))$ \hfill by  (\ref{paso13}) and (\ref{paso28})  \label{paso99}

\item $x \wedge (((x \too 0) \too (0 \wedge 0)) \too 0) = ((x \too 0) \too 0) \wedge (x \too (x \too 0))$ \hfill by  (\ref{paso30}) and (\ref{paso99})  \label{paso100}

\item $x \wedge (((x \too 0) \too 0) \too 0) = ((x \too 0) \too 0) \wedge (x \too (x \too 0))$ \hfill by  (\ref{paso24}) and (\ref{paso100})  \label{paso101}

\item $x \wedge (x \too 0) = ((x \too 0) \too 0) \wedge (x \too (x \too 0))$ \hfill by  (\ref{paso30}) and (\ref{paso101})  \label{paso102}

\item $x \wedge (x \too 0) = x \wedge (x \too (x \too 0))$ \hfill by  (\ref{paso30}) and (\ref{paso102})  \label{paso103}

\item $(x \too 0) \wedge (x \too y) = x \too ((0 \wedge (x \too 0)) \wedge y)$ \hfill by  (\ref{paso51}) and (\ref{paso9})  \label{paso111}

\item $x \too (0 \wedge y) = x \too ((0 \wedge (x \too 0)) \wedge y)$ \hfill by  (\ref{paso9}) and (\ref{paso111})  \label{paso112}

\item $(x \too x) \wedge (0 \too y) = 0 \too ((x \wedge (x \too 0)) \wedge y)$ \hfill by  (\ref{paso72}) and (\ref{paso9})  \label{paso118}

\item $(x \too x) \wedge y = y \wedge (z \too z)$ \hfill by  (\ref{paso22}) and (\ref{paso23}) \label{paso120}

\item $(x \too x) \wedge y = (y \too 0) \too ((y \too 0) \wedge 0)$ \hfill by  (\ref{paso30}) and (\ref{paso65})  \label{paso121}

\item $(x \too x) \wedge y = (y \too 0) \too (0 \wedge (y \too 0))$ \hfill by  (\ref{paso23}) and (\ref{paso121})  \label{paso122}

\item $(x \too x) \wedge (y \too (z \too z)) = 0 \too (0 \wedge (y \too 0))$ \hfill by  (\ref{paso31}) and (\ref{paso65})  \label{paso123}

\item $((x \too (y \too y)) \too 0) \too (0 \wedge ((x \too (y \too y)) \too 0)) = 0 \too (0 \wedge (x \too 0))$ \hfill by  (\ref{paso122}) and (\ref{paso123})  \label{paso124}

\item $(x \too x) \wedge ((y \too y) \too z) = (y \too y) \too (z \wedge (u \too u))$ \hfill by  (\ref{paso120}) and (\ref{paso65})  \label{paso125}

\item $(x \too x) \wedge z = (y \too y) \too (z \wedge (u \too u))$ \hfill by  (\ref{paso4}) and (\ref{paso125})  \label{paso126}

\item $(z \too 0) \too (0 \wedge (z \too 0)) = (y \too y) \too (z \wedge (u \too u))$ \hfill by  (\ref{paso122}) and (\ref{paso126})  \label{paso127}

\item $(z \too 0) \too (0 \wedge (z \too 0)) = z \wedge (u \too u)$ \hfill by  (\ref{paso4}) and (\ref{paso127})  \label{paso128}

\item $x \wedge (y \too y) = (x \too 0) \too (0 \wedge (x \too 0))$ \hfill by (\ref{paso128})  \label{paso129}

\item $(x \too x) \wedge (y \too (z \too z)) = y \too ((u \too u) \wedge y)$ \hfill by  (\ref{paso120}) and (\ref{paso65})  \label{paso130}

\item $((y \too (z \too z)) \too 0) \too (0 \wedge ((y \too (z \too z)) \too 0)) = y \too ((u \too u) \wedge y)$ \hfill by  (\ref{paso122}) and (\ref{paso130})  \label{paso131}

\item $0 \too (0 \wedge (y \too 0)) = y \too ((u \too u) \wedge y)$ \hfill by  (\ref{paso124}) and (\ref{paso131})  \label{paso132}

\item $0 \too (0 \wedge (y \too 0)) = y \too ((y \too 0) \too (0 \wedge (y \too 0)))$ \hfill by  (\ref{paso122}) and (\ref{paso132})  \label{paso133}

\item $0 \too ((x \wedge (x \too 0)) \wedge y) = ((0 \too y) \too 0) \too (0 \wedge ((0 \too y) \too 0))$ \hfill by  (\ref{paso122}) and (\ref{paso118})  \label{paso135}

\item $(x \wedge y) \wedge (x \too (x \too (y \too (y \too 0)))) = (x \wedge y) \wedge ((x \wedge y) \too 0)$ \hfill by  (\ref{paso14}) and (\ref{paso103})  \label{paso136}

\item $(x \too (x \too (y \too (y \too 0)))) \wedge (x \wedge y) = (x \wedge y) \wedge ((x \wedge y) \too 0)$ \hfill by  (\ref{paso23}) and (\ref{paso136})  \label{paso137}

\item $(x \too x) \wedge (((y \too 0) \wedge (x \too x)) \too 0) = (x \too x) \wedge y$ \hfill by  (\ref{paso4}) and (\ref{paso73})  \label{paso139}

\item $(x \too x) \wedge ((y \too (0 \wedge y)) \too 0) = (x \too x) \wedge y$ \hfill by  (\ref{paso66}) and (\ref{paso139})  \label{paso140}

\item $(((y \too (0 \wedge y)) \too 0) \too 0) \too (0 \wedge (((y \too (0 \wedge y)) \too 0) \too 0)) = (x \too x) \wedge y$ \hfill by  (\ref{paso122}) and (\ref{paso140})  \label{paso141}

\item $((0 \too 0) \too (y \too (0 \wedge y))) \too (0 \wedge (((y \too (0 \wedge y)) \too 0) \too 0)) = (x \too x) \wedge y$ \hfill by  (\ref{paso90}) and (\ref{paso141})  \label{paso142}

\item $(y \too (0 \wedge y)) \too (0 \wedge (((y \too (0 \wedge y)) \too 0) \too 0)) = (x \too x) \wedge y$ \hfill by  (\ref{paso4}) and (\ref{paso142})  \label{paso143}

\item $(y \too (0 \wedge y)) \too (0 \wedge ((0 \too 0) \too (y \too (0 \wedge y)))) = (x \too x) \wedge y$ \hfill by  (\ref{paso90}) and (\ref{paso143})  \label{paso144}

\item $(y \too (0 \wedge y)) \too (0 \wedge (y \too (0 \wedge y))) = (x \too x) \wedge y$ \hfill by  (\ref{paso4}) and (\ref{paso144})  \label{paso145}

\item $(x \too (0 \wedge x)) \too (0 \wedge (x \too (0 \wedge x))) = (x \too 0) \too (0 \wedge (x \too 0))$ \hfill by  (\ref{paso122}) and (\ref{paso145})  \label{paso146}

\item $x \wedge ((y \wedge x) \too 0) = x \wedge (x \too (y \too 0))$ \hfill by  (\ref{paso30}) and (\ref{paso73})  \label{paso147}

\item $x \too (x \too (0 \too (x \too 0))) = x \too (x \too x)$ \hfill by  (\ref{paso31}) and (\ref{paso79})  \label{paso153}

\item $(x \too (y \wedge z)) \too (u \too u) = (x \too y) \too ((x \too y) \too ((x \too z) \too ((x \too z) \too (x \too (y \wedge z)))))$ \\\null\hfill by  (\ref{paso22}) and (\ref{paso38})  \label{paso154}

\item $0 \too ((x \too (y \wedge z)) \too 0) = (x \too y) \too ((x \too y) \too ((x \too z) \too ((x \too z) \too (x \too (y \wedge z)))))$ \\\null\hfill by  (\ref{paso31}) and (\ref{paso154})  \label{paso155}

\item $(x \too (y \wedge z)) \too ((x \too y) \too ((x \too y) \too ((x \too z) \too ((x \too z) \too u)))) = (x \too y) \too ((x \too y) \too ((x \too z) \too ((x \too z) \too ((x \too (y \wedge z)) \too u))))$ \hfill by  (\ref{paso38}) and (\ref{paso38})  \label{paso156}

\item $(0 \wedge x) \wedge ((0 \wedge x) \wedge (((0 \wedge x) \too 0) \too x)) = 0 \wedge x$ \hfill by  (\ref{paso86}) and (\ref{paso98})  \label{paso157}

\item $(0 \wedge x) \wedge ((0 \wedge x) \wedge ((x \too 0) \too (0 \wedge x))) = 0 \wedge x$ \hfill by  (\ref{paso90}) and (\ref{paso157})  \label{paso158}

\item $(0 \wedge x) \wedge ((0 \wedge x) \wedge x) = 0 \wedge x$ \hfill by  (\ref{paso95}) and (\ref{paso158})  \label{paso159}

\item $(0 \wedge x) \wedge (x \wedge (x \wedge 0)) = 0 \wedge x$ \hfill by  (\ref{paso23}) and (\ref{paso159})  \label{paso162}

\item $x \wedge ((x \too 0) \too (0 \too (y \too y))) = (x \too 0) \too (0 \wedge (0 \too 0))$ \hfill by  (\ref{paso83}) and (\ref{paso86})  \label{paso161}

\item $(0 \wedge x) \too (((0 \wedge x) \wedge (x \wedge (x \wedge 0))) \too y) = (0 \wedge x) \too ((0 \wedge x) \too ((x \wedge (x \wedge 0)) \too ((x \wedge (x \wedge 0)) \too y)))$\\\null \hfill by  (\ref{paso162}) and (\ref{paso14})  \label{paso163}

\item $(0 \wedge x) \too ((0 \wedge x) \too y) = (0 \wedge x) \too ((0 \wedge x) \too ((x \wedge (x \wedge 0)) \too ((x \wedge (x \wedge 0)) \too y)))$ \\\null \hfill by  (\ref{paso162}) and (\ref{paso163})  \label{paso164}

\item $0 \too (0 \too (x \too (x \too y))) = (0 \wedge x) \too ((0 \wedge x) \too ((x \wedge (x \wedge 0)) \too ((x \wedge (x \wedge 0)) \too y)))$ \\\null \hfill by  (\ref{paso14}) and (\ref{paso164})  \label{paso165}

\item $0 \too (0 \too (x \too (x \too y))) = (0 \wedge x) \too ((0 \wedge x) \too (x \too (x \too ((x \wedge 0) \too ((x \wedge 0) \too y)))))$ \\\null\hfill by  (\ref{paso14}) and (\ref{paso165})  \label{paso166}

\item $0 \too (0 \too (x \too (x \too y))) = (0 \wedge x) \too ((0 \wedge x) \too (x \too (x \too (x \too (x \too (0 \too (0 \too y)))))))$ \\\null\hfill by  (\ref{paso14}) and (\ref{paso166})  \label{paso167}

\item $0 \too (0 \too (x \too (x \too y))) = (0 \wedge x) \too ((0 \wedge x) \too (x \too (x \too (0 \too (0 \too y)))))$ \\\null\hfill by  (\ref{paso79}) and (\ref{paso167})  \label{paso168}

\item $0 \too (0 \too (x \too (x \too y))) = 0 \too (0 \too (x \too (x \too (x \too (x \too (0 \too (0 \too y)))))))$ \\\null\hfill by  (\ref{paso14}) and (\ref{paso168})  \label{paso169}

\item $0 \too (0 \too (x \too (x \too y))) = 0 \too (0 \too (x \too (x \too (0 \too (0 \too y)))))$ \hfill by  (\ref{paso79}) and (\ref{paso169})  \label{paso170}

\item $(0 \wedge x) \too ((0 \wedge ((0 \wedge x) \too 0)) \wedge x) = y \too y$ \hfill by  (\ref{paso22}) and (\ref{paso112})  \label{paso172}

\item $(0 \wedge x) \too ((0 \wedge (0 \too (x \too 0))) \wedge x) = y \too y$ \hfill by  (\ref{paso85}) and (\ref{paso172})  \label{paso173}

\item $(0 \wedge x) \too (x \wedge (0 \wedge (0 \too (x \too 0)))) = y \too y$ \hfill by  (\ref{paso23}) and (\ref{paso173})  \label{paso174}

\item $x \wedge ((x \too 0) \too (0 \too (y \too y))) = (x \too 0) \too 0$ \hfill by  (\ref{paso182}) and (\ref{paso161})  \label{paso184}

\item $x \wedge ((x \too 0) \too (0 \too (y \too y))) = x$ \hfill by  (\ref{paso30}) and (\ref{paso184})  \label{paso185}

\item $0 \wedge (0 \too (x \too x)) = 0$ \hfill by  (\ref{paso182}) and (\ref{paso83})  \label{paso186}

\item $((x \too 0) \too 0) \wedge (0 \too x) = (x \too 0) \too 0$ \hfill by  (\ref{paso182}) and (\ref{paso36})  \label{paso187}

\item $x \wedge (0 \too x) = (x \too 0) \too 0$ \hfill by  (\ref{paso30}) and (\ref{paso187})  \label{paso188}

\item $x \wedge (0 \too x) = x$ \hfill by  (\ref{paso30}) and (\ref{paso188})  \label{paso189}

\item $((x \too (0 \too (x \too 0))) \too 0) \too (x \too 0) = x \too (0 \too (x \too 0))$ \hfill by  (\ref{paso189}) and (\ref{paso54})  \label{paso190}

\item $x \too (((x \too (0 \too (x \too 0))) \too 0) \too 0) = x \too (0 \too (x \too 0))$ \hfill by  (\ref{paso89}) and (\ref{paso190})  \label{paso191}

\item $x \too ((0 \too 0) \too (x \too (0 \too (x \too 0)))) = x \too (0 \too (x \too 0))$ \hfill by  (\ref{paso90}) and (\ref{paso191})  \label{paso192}

\item $x \too (x \too (0 \too (x \too 0))) = x \too (0 \too (x \too 0))$ \hfill by  (\ref{paso4}) and (\ref{paso192})  \label{paso193}

\item $x \too (x \too x) = x \too (0 \too (x \too 0))$ \hfill by  (\ref{paso153}) and (\ref{paso193})  \label{paso194}

\item $(0 \too (0 \too (0 \too x))) \wedge (0 \too (0 \too x)) = 0 \too (0 \too (0 \too x))$ \hfill by  (\ref{paso79}) and (\ref{paso189})  \label{paso196}

\item $(0 \too (0 \too x)) \wedge (0 \too (0 \too (0 \too x))) = 0 \too (0 \too (0 \too x))$ \hfill by  (\ref{paso23}) and (\ref{paso196})  \label{paso197}

\item $0 \too (0 \too x) = 0 \too (0 \too (0 \too x))$ \hfill by  (\ref{paso189}) and (\ref{paso197})  \label{paso198}

\item $x \wedge (((0 \too (y \too y)) \too 0) \too ((x \too 0) \too 0)) = x$ \hfill by  (\ref{paso13}) and (\ref{paso185})  \label{paso200}

\item $x \wedge (((0 \too (y \too y)) \too 0) \too x) = x$ \hfill by  (\ref{paso30}) and (\ref{paso200})  \label{paso201}

\item $((x \too x) \too 0) \too (x \too 0) = x \too (0 \too (x \too 0))$ \hfill by  (\ref{paso194}) and (\ref{paso13})  \label{paso202}

\item $0 \too (x \too 0) = x \too (0 \too (x \too 0))$ \hfill by  (\ref{paso4}) and (\ref{paso202})  \label{paso203}

\item $(x \too 0) \too (0 \too ((x \too 0) \too 0)) = (x \too 0) \too (x \too x)$ \hfill by  (\ref{paso13}) and (\ref{paso194})  \label{paso205}

\item $(x \too 0) \too (0 \too x) = (x \too 0) \too (x \too x)$ \hfill by  (\ref{paso30}) and (\ref{paso205})  \label{paso206}

\item $(x \too 0) \too (0 \too x) = 0 \too ((x \too 0) \too 0)$ \hfill by  (\ref{paso31}) and (\ref{paso206})  \label{paso207}

\item $(x \too 0) \too (0 \too x) = 0 \too x$ \hfill by  (\ref{paso30}) and (\ref{paso207})  \label{paso208}

\item $(x \wedge y) \too ((x \wedge y) \too ((x \wedge y) \too (x \wedge y))) = x \too (x \too (y \too (y \too (0 \too ((x \wedge y) \too 0)))))$ \\\null\hfill by  (\ref{paso194}) and (\ref{paso14})  \label{paso209}

\item $(x \wedge y) \too (0 \too ((x \wedge y) \too 0)) = x \too (x \too (y \too (y \too (0 \too ((x \wedge y) \too 0)))))$ \\\null\hfill by  (\ref{paso31}) and (\ref{paso209})  \label{paso210}

\item $(x \wedge y) \too (x \too (x \too (y \too (y \too (x \wedge y))))) = x \too (x \too (y \too (y \too (0 \too ((x \wedge y) \too 0)))))$\\\null\hfill by  (\ref{paso68}) and (\ref{paso210})  \label{paso211}

\item $x \too (x \too (y \too (y \too ((x \wedge y) \too (x \wedge y))))) = x \too (x \too (y \too (y \too (0 \too ((x \wedge y) \too 0)))))$ \\\null\hfill by  (\ref{paso39}) and (\ref{paso211})  \label{paso212}

\item $x \too (x \too (y \too (0 \too (y \too 0)))) = x \too (x \too (y \too (y \too (0 \too ((x \wedge y) \too 0)))))$\\\null \hfill by  (\ref{paso31}) and (\ref{paso212})  \label{paso213}

\item $x \too (x \too (0 \too (y \too 0))) = x \too (x \too (y \too (y \too (0 \too ((x \wedge y) \too 0)))))$ \hfill by  (\ref{paso203}) and (\ref{paso213})  \label{paso214}

\item $x \too (x \too (0 \too (y \too 0))) = x \too (x \too (y \too (y \too (x \too (x \too (y \too (y \too (x \wedge y))))))))$ \\\null\hfill by  (\ref{paso68}) and (\ref{paso214})  \label{paso215}

\item $x \too (x \too (0 \too (y \too 0))) = x \too (x \too (y \too (y \too (x \wedge y))))$ \hfill by  (\ref{paso45}) and (\ref{paso215})  \label{paso216}

\item $(x \too (y \wedge z)) \too ((x \too (y \wedge z)) \too ((x \too (y \wedge z)) \too (x \too (y \wedge z)))) = (x \too y) \too ((x \too y) \too ((x \too z) \too ((x \too z) \too (0 \too ((x \too (y \wedge z)) \too 0)))))$ \hfill by  (\ref{paso194}) and (\ref{paso38})  \label{paso218}

\item $(x \too (y \wedge z)) \too (0 \too ((x \too (y \wedge z)) \too 0)) = (x \too y) \too ((x \too y) \too ((x \too z) \too ((x \too z) \too (0 \too ((x \too (y \wedge z)) \too 0)))))$ \hfill by  (\ref{paso31}) and (\ref{paso218})  \label{paso219}

\item $(x \too (y \wedge z)) \too ((x \too y) \too ((x \too y) \too ((x \too z) \too ((x \too z) \too (x \too (y \wedge z)))))) = (x \too y) \too ((x \too y) \too ((x \too z) \too ((x \too z) \too (0 \too ((x \too (y \wedge z)) \too 0)))))$ \hfill by  (\ref{paso155}) and (\ref{paso219})  \label{paso220}

\item $(x \too y) \too ((x \too y) \too ((x \too z) \too ((x \too z) \too ((x \too (y \wedge z)) \too (x \too (y \wedge z)))))) = (x \too y) \too ((x \too y) \too ((x \too z) \too ((x \too z) \too (0 \too ((x \too (y \wedge z)) \too 0)))))$ \hfill by  (\ref{paso156}) and (\ref{paso220})  \label{paso221}

\item $(x \too y) \too ((x \too y) \too ((x \too z) \too (0 \too ((x \too z) \too 0)))) = (x \too y) \too ((x \too y) \too ((x \too z) \too ((x \too z) \too (0 \too ((x \too (y \wedge z)) \too 0)))))$ \hfill by  (\ref{paso31}) and (\ref{paso221})  \label{paso222}

\item $(x \too y) \too ((x \too y) \too (0 \too ((x \too z) \too 0))) = (x \too y) \too ((x \too y) \too ((x \too z) \too ((x \too z) \too (0 \too ((x \too (y \wedge z)) \too 0)))))$ \hfill by  (\ref{paso203}) and (\ref{paso222})  \label{paso223}

\item $(x \too y) \too ((x \too y) \too (0 \too ((x \too z) \too 0))) = (x \too y) \too ((x \too y) \too ((x \too z) \too ((x \too z) \too ((x \too y) \too ((x \too y) \too ((x \too z) \too ((x \too z) \too (x \too (y \wedge z)))))))))$ \\\null\hfill by  (\ref{paso155}) and (\ref{paso223})  \label{paso224}

\item $(x \too y) \too ((x \too y) \too (0 \too ((x \too z) \too 0))) = (x \too y) \too ((x \too y) \too ((x \too z) \too ((x \too z) \too (x \too (y \wedge z)))))$ \hfill by  (\ref{paso45}) and (\ref{paso224})  \label{paso225}

\item $0 \too ((x \wedge y) \too 0) = x \too (x \too (0 \too (y \too 0)))$ \hfill by  (\ref{paso216}) and (\ref{paso68})  \label{paso227}

\item $0 \too ((x \too (y \wedge z)) \too 0) = (x \too y) \too ((x \too y) \too (0 \too ((x \too z) \too 0)))$ \hfill by  (\ref{paso225}) and (\ref{paso155})  \label{paso228}

\item $(0 \too x) \wedge (((0 \too x) \too 0) \too y) = ((0 \too x) \too 0) \too (((x \too 0) \too 0) \wedge y)$ \hfill by  (\ref{paso208}) and (\ref{paso35})  \label{paso229}

\item $(0 \too x) \wedge (((0 \too x) \too 0) \too y) = ((0 \too x) \too 0) \too (x \wedge y)$ \hfill by  (\ref{paso30}) and (\ref{paso229})  \label{paso230}

\item $(0 \too (0 \too x)) \wedge (0 \too y) = 0 \too ((0 \too (0 \too x)) \wedge y)$ \hfill by  (\ref{paso198}) and (\ref{paso9})  \label{paso231}

\item $0 \too ((0 \too x) \wedge y) = 0 \too ((0 \too (0 \too x)) \wedge y)$ \hfill by  (\ref{paso9}) and (\ref{paso231})  \label{paso232}

\item $(x \wedge 0) \wedge ((0 \too (y \too y)) \wedge (((0 \too (y \too y)) \too 0) \too x)) = x \wedge 0$ \hfill by  (\ref{paso88}) and (\ref{paso201})  \label{paso237}

\item $(x \wedge 0) \wedge (((0 \too (y \too y)) \too 0) \too ((y \too y) \wedge x)) = x \wedge 0$ \hfill by  (\ref{paso230}) and (\ref{paso237})  \label{paso238}

\item $(x \wedge 0) \wedge (((0 \too (y \too y)) \too 0) \too ((x \too 0) \too (0 \wedge (x \too 0)))) = x \wedge 0$ \hfill by  (\ref{paso122}) and (\ref{paso238})  \label{paso239}

\item $((0 \wedge (0 \too 0)) \too 0) \wedge ((x \wedge ((0 \too 0) \too 0)) \too 0) = (((0 \too 0) \too 0) \wedge (((x \too 0) \wedge ((0 \too 0) \too 0)) \too 0)) \too 0$ \hfill by  (\ref{paso76}) and (\ref{paso57})  \label{paso240}

\item $(0 \too 0) \wedge ((x \wedge ((0 \too 0) \too 0)) \too 0) = (((0 \too 0) \too 0) \wedge (((x \too 0) \wedge ((0 \too 0) \too 0)) \too 0)) \too 0$ \\\null\hfill by  (\ref{paso189}) and (\ref{paso240})  \label{paso241}

\item $(0 \too 0) \wedge ((x \wedge 0) \too 0) = (((0 \too 0) \too 0) \wedge (((x \too 0) \wedge ((0 \too 0) \too 0)) \too 0)) \too 0$ \\\null \hfill by  (\ref{paso4}) and (\ref{paso241})  \label{paso242}

\item $(((x \wedge 0) \too 0) \too 0) \too (0 \wedge (((x \wedge 0) \too 0) \too 0)) = (((0 \too 0) \too 0) \wedge (((x \too 0) \wedge ((0 \too 0) \too 0)) \too 0)) \too 0$ \hfill by  (\ref{paso122}) and (\ref{paso242})  \label{paso243}

\item $((0 \too 0) \too (x \wedge 0)) \too (0 \wedge (((x \wedge 0) \too 0) \too 0)) = (((0 \too 0) \too 0) \wedge (((x \too 0) \wedge ((0 \too 0) \too 0)) \too 0)) \too 0$ \hfill by  (\ref{paso90}) and (\ref{paso243})  \label{paso244}

\item $(x \wedge 0) \too (0 \wedge (((x \wedge 0) \too 0) \too 0)) = (((0 \too 0) \too 0) \wedge (((x \too 0) \wedge ((0 \too 0) \too 0)) \too 0)) \too 0$ \\\null \hfill by  (\ref{paso4}) and (\ref{paso244})  \label{paso245}

\item $(x \wedge 0) \too (0 \wedge ((0 \too 0) \too (x \wedge 0))) = (((0 \too 0) \too 0) \wedge (((x \too 0) \wedge ((0 \too 0) \too 0)) \too 0)) \too 0$ \\\null\hfill by  (\ref{paso90}) and (\ref{paso245})  \label{paso246}

\item $(x \wedge 0) \too (0 \wedge (x \wedge 0)) = (((0 \too 0) \too 0) \wedge (((x \too 0) \wedge ((0 \too 0) \too 0)) \too 0)) \too 0$ \\\null\hfill by  (\ref{paso4}) and (\ref{paso246})  \label{paso247}

\item $(x \wedge 0) \too (0 \wedge (x \wedge 0)) = (0 \wedge (((x \too 0) \wedge ((0 \too 0) \too 0)) \too 0)) \too 0$ \hfill by  (\ref{paso4}) and (\ref{paso247})  \label{paso248}

\item $(x \wedge 0) \too (0 \wedge (x \wedge 0)) = (0 \wedge (((x \too 0) \wedge 0) \too 0)) \too 0$ \hfill by  (\ref{paso4}) and (\ref{paso248})  \label{paso249}

\item $(x \wedge 0) \too (0 \wedge (x \wedge 0)) = (0 \wedge ((0 \wedge (x \too 0)) \too 0)) \too 0$ \hfill by  (\ref{paso23}) and (\ref{paso249})  \label{paso250}

\item $(x \wedge 0) \too (0 \wedge (x \wedge 0)) = (0 \wedge (0 \too ((x \too 0) \too 0))) \too 0$ \hfill by  (\ref{paso85}) and (\ref{paso250})  \label{paso251}

\item $(x \wedge 0) \too (0 \wedge (x \wedge 0)) = (0 \wedge (0 \too x)) \too 0$ \hfill by  (\ref{paso30}) and (\ref{paso251})  \label{paso252}

\item $(0 \wedge x) \too (x \wedge 0) = y \too y$ \hfill by  (\ref{paso22}) and (\ref{paso23})  \label{paso254}

\item $(0 \wedge x) \too (y \too y) = 0 \too (0 \too (x \too (x \too (x \wedge 0))))$ \hfill by  (\ref{paso254}) and (\ref{paso14})  \label{paso255}

\item $0 \too ((0 \wedge x) \too 0) = 0 \too (0 \too (x \too (x \too (x \wedge 0))))$ \hfill by  (\ref{paso31}) and (\ref{paso255})  \label{paso256}

\item $0 \too (0 \too (0 \too (x \too 0))) = 0 \too (0 \too (x \too (x \too (x \wedge 0))))$ \hfill by  (\ref{paso227}) and (\ref{paso256})  \label{paso257}

\item $0 \too (0 \too (x \too 0)) = 0 \too (0 \too (x \too (x \too (x \wedge 0))))$ \hfill by  (\ref{paso198}) and (\ref{paso257})  \label{paso258}

\item $(0 \wedge x) \wedge (((x \wedge 0) \too 0) \too ((0 \wedge x) \too 0)) = (0 \wedge x) \wedge (y \too y)$ \hfill by  (\ref{paso254}) and (\ref{paso33})  \label{paso260}

\item $((x \too x) \too 0) \too (((0 \wedge y) \too 0) \wedge (y \wedge 0)) = (0 \wedge y) \too (y \wedge 0)$ \hfill by  (\ref{paso254}) and (\ref{paso54})  \label{paso262}

\item $0 \too (((0 \wedge x) \too 0) \wedge (x \wedge 0)) = (0 \wedge x) \too (x \wedge 0)$ \hfill by  (\ref{paso4}) and (\ref{paso262})  \label{paso263}

\item $(0 \wedge x) \too ((0 \wedge x) \too ((0 \wedge x) \too (y \too y))) = (0 \wedge x) \too ((0 \wedge x) \too (x \wedge 0))$ \hfill by  (\ref{paso254}) and (\ref{paso79})  \label{paso264}

\item $(0 \wedge x) \too ((0 \wedge x) \too (0 \too ((0 \wedge x) \too 0))) = (0 \wedge x) \too ((0 \wedge x) \too (x \wedge 0))$ \hfill by  (\ref{paso31}) and (\ref{paso264})  \label{paso265}

\item $(0 \wedge x) \too ((0 \wedge x) \too (0 \too (0 \too (0 \too (x \too 0))))) = (0 \wedge x) \too ((0 \wedge x) \too (x \wedge 0))$ \\\null\hfill by  (\ref{paso227}) and (\ref{paso265})  \label{paso266}

\item $(0 \wedge x) \too ((0 \wedge x) \too (0 \too (0 \too (x \too 0)))) = (0 \wedge x) \too ((0 \wedge x) \too (x \wedge 0))$ \hfill by  (\ref{paso198}) and (\ref{paso266})  \label{paso267}

\item $0 \too (0 \too (x \too (x \too (0 \too (0 \too (x \too 0)))))) = (0 \wedge x) \too ((0 \wedge x) \too (x \wedge 0))$ \\\null\hfill by  (\ref{paso14}) and (\ref{paso267})  \label{paso268}

\item $0 \too (0 \too (x \too (x \too (x \too 0)))) = (0 \wedge x) \too ((0 \wedge x) \too (x \wedge 0))$ \hfill by  (\ref{paso170}) and (\ref{paso268})  \label{paso269}

\item $0 \too (0 \too (x \too (x \too 0))) = (0 \wedge x) \too ((0 \wedge x) \too (x \wedge 0))$ \hfill by  (\ref{paso82}) and (\ref{paso269})  \label{paso270}

\item $0 \too (0 \too (x \too (x \too 0))) = 0 \too (0 \too (x \too (x \too (x \wedge 0))))$ \hfill by  (\ref{paso14}) and (\ref{paso270})  \label{paso271}

\item $0 \too (0 \too (x \too (x \too 0))) = 0 \too (0 \too (x \too 0))$ \hfill by  (\ref{paso258}) and (\ref{paso271})  \label{paso272}

\item $((x \wedge y) \too 0) \wedge (y \wedge x) = (x \wedge y) \wedge ((x \wedge y) \too 0)$ \hfill by (\ref{paso23}) and (\ref{paso23}) \label{paso276}

\item $((x \wedge y) \too 0) \wedge (y \wedge x) = (x \too (x \too (y \too (y \too 0)))) \wedge (x \wedge y)$ \hfill by  (\ref{paso137}) and (\ref{paso276})  \label{paso277}

\item $0 \too ((0 \too (0 \too (x \too (x \too 0)))) \wedge (0 \wedge x)) = (0 \wedge x) \too (x \wedge 0)$ \hfill by  (\ref{paso277}) and (\ref{paso263})  \label{paso278}

\item $0 \too ((0 \too (0 \too (x \too 0))) \wedge (0 \wedge x)) = (0 \wedge x) \too (x \wedge 0)$ \hfill by  (\ref{paso272}) and (\ref{paso278})  \label{paso279}

\item $0 \too ((0 \too (x \too 0)) \wedge (0 \wedge x)) = (0 \wedge x) \too (x \wedge 0)$ \hfill by  (\ref{paso232}) and (\ref{paso279})  \label{paso280}

\item $(0 \too 0) \wedge (((x \too 0) \wedge (0 \too 0)) \too 0) = (((0 \wedge (0 \too (0 \too 0))) \too 0) \wedge ((x \wedge ((0 \wedge ((0 \too 0) \too 0)) \too 0)) \too 0)) \too 0$ \hfill by  (\ref{paso98}) and (\ref{paso60})  \label{paso281}

\item $(0 \too 0) \wedge ((x \too (0 \wedge x)) \too 0) = (((0 \wedge (0 \too (0 \too 0))) \too 0) \wedge ((x \wedge ((0 \wedge ((0 \too 0) \too 0)) \too 0)) \too 0)) \too 0$ \hfill by  (\ref{paso66}) and (\ref{paso281})  \label{paso282}

\item $(((x \too (0 \wedge x)) \too 0) \too 0) \too (0 \wedge (((x \too (0 \wedge x)) \too 0) \too 0)) = (((0 \wedge (0 \too (0 \too 0))) \too 0) \wedge ((x \wedge ((0 \wedge ((0 \too 0) \too 0)) \too 0)) \too 0)) \too 0$ \hfill by  (\ref{paso122}) and (\ref{paso282})  \label{paso283}

\item $((0 \too 0) \too (x \too (0 \wedge x))) \too (0 \wedge (((x \too (0 \wedge x)) \too 0) \too 0)) = (((0 \wedge (0 \too (0 \too 0))) \too 0) \wedge ((x \wedge ((0 \wedge ((0 \too 0) \too 0)) \too 0)) \too 0)) \too 0$ \hfill by  (\ref{paso90}) and (\ref{paso283})  \label{paso284}

\item $(x \too (0 \wedge x)) \too (0 \wedge (((x \too (0 \wedge x)) \too 0) \too 0)) = (((0 \wedge (0 \too (0 \too 0))) \too 0) \wedge ((x \wedge ((0 \wedge ((0 \too 0) \too 0)) \too 0)) \too 0)) \too 0$ \hfill by  (\ref{paso4}) and (\ref{paso284})  \label{paso285}

\item $(x \too (0 \wedge x)) \too (0 \wedge ((0 \too 0) \too (x \too (0 \wedge x)))) = (((0 \wedge (0 \too (0 \too 0))) \too 0) \wedge ((x \wedge ((0 \wedge ((0 \too 0) \too 0)) \too 0)) \too 0)) \too 0$ \hfill by  (\ref{paso90}) and (\ref{paso285})  \label{paso286}

\item $(x \too (0 \wedge x)) \too (0 \wedge (x \too (0 \wedge x))) = (((0 \wedge (0 \too (0 \too 0))) \too 0) \wedge ((x \wedge ((0 \wedge ((0 \too 0) \too 0)) \too 0)) \too 0)) \too 0$ \hfill by  (\ref{paso4}) and (\ref{paso286})  \label{paso287}

\item $(x \too 0) \too (0 \wedge (x \too 0)) = (((0 \wedge (0 \too (0 \too 0))) \too 0) \wedge ((x \wedge ((0 \wedge ((0 \too 0) \too 0)) \too 0)) \too 0)) \too 0$ \hfill by  (\ref{paso146}) and (\ref{paso287})  \label{paso288}

\item $(x \too 0) \too (0 \wedge (x \too 0)) = (((0 \wedge (0 \too 0)) \too 0) \wedge ((x \wedge ((0 \wedge ((0 \too 0) \too 0)) \too 0)) \too 0)) \too 0$ \\\null\hfill by  (\ref{paso103}) and (\ref{paso288})  \label{paso289}

\item $(x \too 0) \too (0 \wedge (x \too 0)) = ((0 \too 0) \wedge ((x \wedge ((0 \wedge ((0 \too 0) \too 0)) \too 0)) \too 0)) \too 0$ \\\null\hfill by  (\ref{paso189}) and (\ref{paso289})  \label{paso290}

\item $(x \too 0) \too (0 \wedge (x \too 0)) = ((0 \too 0) \wedge ((x \wedge ((0 \wedge 0) \too 0)) \too 0)) \too 0$ \hfill by  (\ref{paso4}) and (\ref{paso290})  \label{paso291}

\item $(x \too 0) \too (0 \wedge (x \too 0)) = ((0 \too 0) \wedge ((x \wedge (0 \too 0)) \too 0)) \too 0$ \hfill by  (\ref{paso24}) and (\ref{paso291})  \label{paso292}

\item $(x \too 0) \too (0 \wedge (x \too 0)) = ((0 \too 0) \wedge ((0 \too 0) \too (x \too 0))) \too 0$ \hfill by  (\ref{paso147}) and (\ref{paso292})  \label{paso293}

\item $(x \too 0) \too (0 \wedge (x \too 0)) = ((0 \too 0) \wedge (x \too 0)) \too 0$ \hfill by  (\ref{paso4}) and (\ref{paso293})  \label{paso294}

\item $(x \too 0) \too (0 \wedge (x \too 0)) = (((x \too 0) \too 0) \too (0 \wedge ((x \too 0) \too 0))) \too 0$ \hfill by  (\ref{paso122}) and (\ref{paso294})  \label{paso295}

\item $(x \too 0) \too (0 \wedge (x \too 0)) = (x \too (0 \wedge ((x \too 0) \too 0))) \too 0$ \hfill by  (\ref{paso30}) and (\ref{paso295})  \label{paso296}

\item $(x \too 0) \too (0 \wedge (x \too 0)) = (x \too (0 \wedge x)) \too 0$ \hfill by  (\ref{paso30}) and (\ref{paso296})  \label{paso297}

\item $(x \wedge 0) \wedge (((0 \too (y \too y)) \too 0) \too ((x \too (0 \wedge x)) \too 0)) = x \wedge 0$ \hfill by  (\ref{paso297}) and (\ref{paso239})  \label{paso298}

\item $0 \too ((x \wedge (x \too 0)) \wedge y) = ((0 \too y) \too (0 \wedge (0 \too y))) \too 0$ \hfill by  (\ref{paso297}) and (\ref{paso135})  \label{paso300}

\item $x \too ((x \too (0 \wedge x)) \too 0) = 0 \too (0 \wedge (x \too 0))$ \hfill by  (\ref{paso297}) and (\ref{paso133})  \label{paso301}

\item $x \wedge (y \too y) = (x \too (0 \wedge x)) \too 0$ \hfill by  (\ref{paso297}) and (\ref{paso129})  \label{paso303}

\item $x \wedge (y \too y) = (x \too (x \wedge 0)) \too 0$ \hfill by  (\ref{paso23}) and (\ref{paso303})  \label{paso305}

\item $x \wedge (((x \too (0 \wedge x)) \too 0) \too 0) = x \wedge (x \too ((y \too y) \too 0))$ \hfill by  (\ref{paso303}) and (\ref{paso85})  \label{paso306}

\item $x \wedge ((0 \too 0) \too (x \too (0 \wedge x))) = x \wedge (x \too ((y \too y) \too 0))$ \hfill by  (\ref{paso90}) and (\ref{paso306})  \label{paso307}

\item $x \wedge (x \too (0 \wedge x)) = x \wedge (x \too ((y \too y) \too 0))$ \hfill by  (\ref{paso4}) and (\ref{paso307})  \label{paso308}

\item $x \wedge (x \too (0 \wedge x)) = x \wedge (x \too 0)$ \hfill by  (\ref{paso4}) and (\ref{paso308})  \label{paso309}

\item $(0 \wedge x) \wedge (((x \wedge 0) \too 0) \too ((0 \wedge x) \too 0)) = ((0 \wedge x) \too ((0 \wedge x) \wedge 0)) \too 0$ \hfill by  (\ref{paso305}) and (\ref{paso260})  \label{paso311}

\item $(0 \wedge x) \wedge (((x \wedge 0) \too 0) \too ((0 \wedge x) \too 0)) = ((0 \wedge x) \too (0 \wedge (0 \wedge x))) \too 0$ \hfill by  (\ref{paso23}) and (\ref{paso311})  \label{paso312}

\item $((((x \too x) \too 0) \too (x \too x)) \too (0 \wedge (((x \too x) \too 0) \too (x \too x)))) \too 0 = x \too x$ \\\null\hfill by  (\ref{paso303}) and (\ref{paso98})  \label{paso313}

\item $((0 \too (x \too x)) \too (0 \wedge (((x \too x) \too 0) \too (x \too x)))) \too 0 = x \too x$ \hfill by  (\ref{paso4}) and (\ref{paso313})  \label{paso314}

\item $((0 \too (x \too x)) \too (0 \wedge (0 \too (x \too x)))) \too 0 = x \too x$ \hfill by  (\ref{paso4}) and (\ref{paso314})  \label{paso315}

\item $((0 \too (x \too x)) \too 0) \too 0 = x \too x$ \hfill by  (\ref{paso186}) and (\ref{paso315})  \label{paso316}

\item $0 \too (x \too x) = x \too x$ \hfill by  (\ref{paso30}) and (\ref{paso316})  \label{paso318}

\item $((0 \wedge x) \too (x \wedge 0)) \wedge y = (y \too (0 \wedge y)) \too 0$ \hfill by  (\ref{paso254}) and (\ref{paso303})  \label{paso319}

\item $(x \wedge 0) \wedge (((y \too y) \too 0) \too ((x \too (0 \wedge x)) \too 0)) = x \wedge 0$ \hfill by  (\ref{paso318}) and (\ref{paso298})  \label{paso328}

\item $(x \wedge 0) \wedge (0 \too ((x \too (0 \wedge x)) \too 0)) = x \wedge 0$ \hfill by  (\ref{paso4}) and (\ref{paso328})  \label{paso329}

\item $(x \wedge 0) \wedge ((x \too 0) \too ((x \too 0) \too (0 \too ((x \too x) \too 0)))) = x \wedge 0$ \hfill by  (\ref{paso228}) and (\ref{paso329})  \label{paso330}

\item $(x \wedge 0) \wedge ((x \too 0) \too ((x \too 0) \too (0 \too 0))) = x \wedge 0$ \hfill by  (\ref{paso4}) and (\ref{paso330})  \label{paso331}

\item $(x \wedge 0) \wedge ((x \too 0) \too (0 \too ((x \too 0) \too 0))) = x \wedge 0$ \hfill by  (\ref{paso31}) and (\ref{paso331})  \label{paso332}

\item $(x \wedge 0) \wedge ((x \too 0) \too (0 \too x)) = x \wedge 0$ \hfill by  (\ref{paso30}) and (\ref{paso332})  \label{paso333}

\item $(x \wedge 0) \wedge (0 \too x) = x \wedge 0$ \hfill by  (\ref{paso208}) and (\ref{paso333})  \label{paso334}

\item $(0 \too x) \wedge (x \wedge 0) = x \wedge 0$ \hfill by  (\ref{paso23}) and (\ref{paso334})  \label{paso335}

\item $(x \too x) \wedge (0 \too y) = 0 \too ((x \too x) \wedge y)$ \hfill by  (\ref{paso318}) and (\ref{paso9})  \label{paso336}

\item $((0 \too y) \too (0 \wedge (0 \too y))) \too 0 = 0 \too ((x \too x) \wedge y)$ \hfill by  (\ref{paso303}) and (\ref{paso336})  \label{paso337}

\item $((0 \too y) \too (0 \wedge (0 \too y))) \too 0 = 0 \too ((y \too (0 \wedge y)) \too 0)$ \hfill by  (\ref{paso303}) and (\ref{paso337})  \label{paso338}

\item $((0 \too y) \too (0 \wedge (0 \too y))) \too 0 = (y \too 0) \too ((y \too 0) \too (0 \too ((y \too y) \too 0)))$\\\null \hfill by  (\ref{paso228}) and (\ref{paso338})  \label{paso339}

\item $((0 \too y) \too (0 \wedge (0 \too y))) \too 0 = (y \too 0) \too ((y \too 0) \too (0 \too 0))$ \hfill by  (\ref{paso4}) and (\ref{paso339})  \label{paso340}

\item $((0 \too y) \too (0 \wedge (0 \too y))) \too 0 = (y \too 0) \too (0 \too ((y \too 0) \too 0))$ \hfill by  (\ref{paso31}) and (\ref{paso340})  \label{paso341}

\item $((0 \too y) \too (0 \wedge (0 \too y))) \too 0 = (y \too 0) \too (0 \too y)$ \hfill by  (\ref{paso30}) and (\ref{paso341})  \label{paso342}

\item $((0 \too x) \too (0 \wedge (0 \too x))) \too 0 = 0 \too x$ \hfill by  (\ref{paso208}) and (\ref{paso342})  \label{paso343}

\item $0 \too ((x \wedge (x \too 0)) \wedge y) = 0 \too y$ \hfill by  (\ref{paso343}) and (\ref{paso300})  \label{paso344}

\item $x \wedge ((x \wedge (x \too 0)) \too 0) = x \wedge (x \too ((x \too (0 \wedge x)) \too 0))$ \hfill by  (\ref{paso309}) and (\ref{paso85})  \label{paso345}

\item $x \wedge (x \too ((x \too 0) \too 0)) = x \wedge (x \too ((x \too (0 \wedge x)) \too 0))$ \hfill by  (\ref{paso85}) and (\ref{paso345})  \label{paso346}

\item $x \wedge (x \too x) = x \wedge (x \too ((x \too (0 \wedge x)) \too 0))$ \hfill by  (\ref{paso30}) and (\ref{paso346})  \label{paso347}

\item $(x \too (x \wedge 0)) \too 0 = x \wedge (x \too ((x \too (0 \wedge x)) \too 0))$ \hfill by  (\ref{paso305}) and (\ref{paso347})  \label{paso348}

\item $(x \wedge ((y \too (0 \wedge 0)) \too 0)) \wedge ((((y \wedge x) \too 0) \wedge ((y \wedge x) \too 0)) \too z) = (((y \wedge x) \too 0) \wedge ((y \wedge x) \too 0)) \too (0 \wedge z)$ \hfill by  (\ref{paso9}) and (\ref{paso61})  \label{paso350}

\item $(x \wedge ((y \too 0) \too 0)) \wedge ((((y \wedge x) \too 0) \wedge ((y \wedge x) \too 0)) \too z) = (((y \wedge x) \too 0) \wedge ((y \wedge x) \too 0)) \too (0 \wedge z)$ \hfill by  (\ref{paso24}) and (\ref{paso350})  \label{paso351}

\item $(x \wedge y) \wedge ((((y \wedge x) \too 0) \wedge ((y \wedge x) \too 0)) \too z) = (((y \wedge x) \too 0) \wedge ((y \wedge x) \too 0)) \too (0 \wedge z)$\\\null \hfill by  (\ref{paso30}) and (\ref{paso351})  \label{paso352}

\item $(x \wedge y) \wedge (((y \wedge x) \too 0) \too z) = (((y \wedge x) \too 0) \wedge ((y \wedge x) \too 0)) \too (0 \wedge z)$ \hfill by  (\ref{paso24}) and (\ref{paso352})  \label{paso353}

\item $(x \wedge y) \wedge (((y \wedge x) \too 0) \too z) = ((y \wedge x) \too 0) \too (0 \wedge z)$ \hfill by  (\ref{paso24}) and (\ref{paso353})  \label{paso354}

\item $((x \wedge 0) \too 0) \too (0 \wedge ((0 \wedge x) \too 0)) = ((0 \wedge x) \too (0 \wedge (0 \wedge x))) \too 0$ \hfill by  (\ref{paso354}) and (\ref{paso312})  \label{paso355}

\item $((x \wedge 0) \too 0) \too (0 \wedge (0 \too (x \too 0))) = ((0 \wedge x) \too (0 \wedge (0 \wedge x))) \too 0$ \hfill by  (\ref{paso85}) and (\ref{paso355})  \label{paso356}

\item $0 \too (((x \too x) \wedge 0) \wedge y) = 0 \too y$ \hfill by  (\ref{paso4}) and (\ref{paso344})  \label{paso358}

\item $0 \too ((0 \wedge (x \too x)) \wedge y) = 0 \too y$ \hfill by  (\ref{paso23}) and (\ref{paso358})  \label{paso359}

\item $0 \too (((0 \too (0 \wedge 0)) \too 0) \wedge y) = 0 \too y$ \hfill by  (\ref{paso305}) and (\ref{paso359})  \label{paso360}

\item $0 \too (((0 \too 0) \too 0) \wedge y) = 0 \too y$ \hfill by  (\ref{paso24}) and (\ref{paso360})  \label{paso361}

\item $0 \too (0 \wedge x) = 0 \too x$ \hfill by  (\ref{paso4}) and (\ref{paso361})  \label{paso362}

\item $0 \too (x \too 0) = x \too ((x \too (0 \wedge x)) \too 0)$ \hfill by  (\ref{paso362}) and (\ref{paso301})  \label{paso363}

\item $x \wedge (0 \too (x \too 0)) = (x \too (x \wedge 0)) \too 0$ \hfill by  (\ref{paso363}) and (\ref{paso348})  \label{paso365}

\item $(0 \too x) \wedge (0 \too y) = 0 \too (x \wedge (0 \wedge y))$ \hfill by  (\ref{paso362}) and (\ref{paso9})  \label{paso366}

\item $0 \too (x \wedge y) = 0 \too (x \wedge (0 \wedge y))$ \hfill by  (\ref{paso9}) and (\ref{paso366})  \label{paso367}

\item $0 \too ((0 \too (x \too 0)) \wedge x) = (0 \wedge x) \too (x \wedge 0)$ \hfill by  (\ref{paso367}) and (\ref{paso280})  \label{paso369}

\item $0 \too (x \wedge (0 \too (x \too 0))) = (0 \wedge x) \too (x \wedge 0)$ \hfill by  (\ref{paso23}) and (\ref{paso369})  \label{paso370}

\item $0 \too ((x \too (x \wedge 0)) \too 0) = (0 \wedge x) \too (x \wedge 0)$ \hfill by  (\ref{paso365}) and (\ref{paso370})  \label{paso371}

\item $(x \too x) \too ((x \too x) \too (0 \too ((x \too 0) \too 0))) = (0 \wedge x) \too (x \wedge 0)$ \hfill by  (\ref{paso228}) and (\ref{paso371})  \label{paso372}

\item $(x \too x) \too ((x \too x) \too (0 \too x)) = (0 \wedge x) \too (x \wedge 0)$ \hfill by  (\ref{paso30}) and (\ref{paso372})  \label{paso373}

\item $(x \too x) \too (0 \too x) = (0 \wedge x) \too (x \wedge 0)$ \hfill by  (\ref{paso4}) and (\ref{paso373})  \label{paso374}

\item $0 \too x = (0 \wedge x) \too (x \wedge 0)$ \hfill by  (\ref{paso4}) and (\ref{paso374})  \label{paso375}

\item $(0 \too x) \wedge y = (y \too (0 \wedge y)) \too 0$ \hfill by  (\ref{paso375}) and (\ref{paso319})  \label{paso377}

\item $((x \wedge 0) \too (0 \wedge (x \wedge 0))) \too 0 = x \wedge 0$ \hfill by  (\ref{paso377}) and (\ref{paso335})  \label{paso379}

\item $((0 \wedge (0 \too x)) \too 0) \too 0 = x \wedge 0$ \hfill by  (\ref{paso252}) and (\ref{paso379})  \label{paso380}

\item $(((0 \too (0 \wedge 0)) \too 0) \too 0) \too 0 = x \wedge 0$ \hfill by  (\ref{paso377}) and (\ref{paso380})  \label{paso381}

\item $(((0 \too 0) \too 0) \too 0) \too 0 = x \wedge 0$ \hfill by  (\ref{paso24}) and (\ref{paso381})  \label{paso382}

\item $(0 \too 0) \too 0 = x \wedge 0$ \hfill by  (\ref{paso4}) and (\ref{paso382})  \label{paso383}

\item $0 = x \wedge 0$. \hfill by  (\ref{paso4}) and (\ref{paso383})  \label{paso384}
\end{enumerate}

\end{document}